\documentclass[reqno,11pt,oneside]{amsart}
\usepackage[utf8]{inputenc}
\usepackage{amsmath}
\usepackage{amsfonts}
\usepackage{amssymb}
\usepackage{graphicx}
\usepackage{verbatim}
\usepackage{mathrsfs}
\usepackage{upref,amsthm,amsxtra,exscale}
\usepackage{cite}
\usepackage{palatino}
\usepackage[left=3cm,  right=3cm,  top=4cm,  bottom=4cm]{geometry}
\usepackage[colorlinks=true,urlcolor=blue,
citecolor=red,linkcolor=blue,linktocpage,pdfpagelabels,
bookmarksnumbered,bookmarksopen]{hyperref}
\usepackage{tikz}
\definecolor{planeL}{RGB}{160,200,220}\definecolor{planeR}{RGB}{40,85,140}
\definecolor{capcol}{RGB}{160,60,60}

\makeatletter
\renewcommand{\paragraph}{\@startsection{paragraph}{4}%
  \z@{1.2ex \@plus .2ex \@minus .2ex}%
  {-1em}%
  {\normalfont\bfseries}}
\makeatother

\newtheorem{theorem}{Theorem}[section]
\newtheorem{corollary}[theorem]{Corollary}
\newtheorem{remark}[theorem]{Remark}

\newtheorem{lemma}[theorem]{Lemma}
\newtheorem{proposition}[theorem]{Proposition}

\numberwithin{equation}{section}

\newcommand{\R}{\mathbb R}
\newcommand{\N}{\mathcal N}
\newcommand{\E}{\mathcal E}
\newcommand{\PP}{\mathcal P}
\newcommand{\dist}{\operatorname{dist}}
\newcommand{\eps}{\varepsilon}

\def\r{\mathbb{R}}
\def\rn{\mathbb{R}^n}
\def\hn{\mathbb{R}^n_+}
\def\h{\mathbb{H}}
\def\z{\mathbb{Z}}
\def\z2{\mathbb{Z}_2}

\def\sn{\mathbb{S}^n}
\def\n{\mathbb{N}}

\def\eps{\varepsilon}
\def\rh{\rightharpoonup}

\def\idm{\int_{\partial M}}

\def\ihn{\int_{\hn}}
\def\idhn{\int_{\partial\hn}}
\def\vp{\varphi}

\def\vr{\varrho}

\def\Ecal{\mathcal{E}}

\def\supp{\mathrm{supp}}
\def\what{\widehat}
\def\tilde{\widetilde}

\def\dist{\mathrm{dist}}

\title[A least-energy nodal solution for the Yamabe problem with negative mean curvature]{A least-energy nodal solution for the Yamabe problem with negative mean curvature on the boundary}
\author{Mónica Clapp}
\address[Mónica Clapp]{Instituto de Matemáticas,
Universidad Nacional Autónoma de México,
Campus Juriquilla,
76230 Querétaro, Qro., Mexico}
\email{monica.clapp@im.unam.mx}

\author{Benedetta Pellacci}
\address[Benedetta Pellacci]{Dipartimento di Matematica e Fisica,
Universit\`a della Campania ``Luigi Vanvitelli'',
Viale Lincoln 5,
81100 Caserta, Italy}
\email{benedetta.pellacci@unicampania.it}

\author{Angela Pistoia}
\address[Angela Pistoia]{Dipartimento SBAI, Sapienza Universit\`a di Roma,
via Antonio Scarpa 16, 00161 Roma, Italy }
\email{angela.pistoia@uniroma1.it}

\date{\today}

\thanks{A. Pistoia  is partially supported by  the MUR-PRIN-20227HX33Z
  ``Pattern formation in nonlinear phenomena'' and  partially by
  INDAM-GNAMPA project ``Problemi di doppia curvatura su variet\`a a
  bordo e legami con le EDP di tipo ellittico''. 
 B. Pellacci is partially supported by  the MUR-PRIN-20227HX33Z
  ``Pattern formation in nonlinear phenomena'' and  partially by
  INDAM-GNAMPA project ``Problemi di ottimizzazione in PDEs da modelli biologici''.}

\subjclass[2010]{53C21, 35J60, 58J60}
\keywords{Escobar problem, Yamabe problem, manifolds with boundary, nodal solutions, conformal geometry, variational methods.
}
\begin{document}

\begin{abstract}
Let $(M,g)$ be a positive compact smooth Riemannian manifold with boundary. We consider the conformally invariant critical problem
\[
\begin{cases}
-\Delta_g u+c_nR_gu=a_n|u|^{\frac4{n-2}}u&\text{in }M,\\
\partial_{\nu_g}u+d_nh_gu=b|u|^{\frac2{n-2}}u&\text{on }\partial M,
\end{cases}
\]
where $R_g$ is the scalar curvature of $M$, $h_g$ is the mean curvature of $\partial M$, and $\nu_g$ is the outward unit normal with respect to the metric $g$, $a_n, c_n$ and $d_n$ are positive constants and $b$ is any real number. A positive solution to this problem gives rise to a metric conformal to $g$ of constant positive scalar curvature in $M$ and constant mean curvature on its boundary. Han and Li showed in \cite{HanLi} that, if $n\geq 7$ and $\partial M$ contains a nonumbilic point, then the problem admits a least-energy positive solution. 
Under the same assumptions, we establish here the existence of a least energy sign-changing solution for $b<0$, complementing our previous result in \cite{CPP} for $b\geq 0$.
\end{abstract}

\maketitle

\section{An overview and the statement of our main result}

Let $(M,g)$ be a compact smooth Riemannian manifold with nonempty boundary, of dimension $n\geq3$. We consider the problem
\begin{equation}\label{P}
\begin{cases}
-\Delta_g u+c_nR_gu=a_n|u|^{p-2}u&\text{on }M,\\[1mm]
\partial_{\nu_g}u+d_nh_gu=b|u|^{q-2}u&\text{on }\partial M,
\end{cases}
\end{equation}
where $R_g$ is the scalar curvature of $M$, $h_g$ is the mean curvature of $\partial M$, and $\nu_g$ is the outward unit normal on $\partial M$ with respect to the metric $g$, and
\[
p:=\frac{2n}{n-2},\qquad q:=\frac{2(n-1)}{n-2},\qquad
 a_n:=n(n-2),\quad c_n:=\frac{n-2}{4(n-1)},\quad d_n:=\frac{n-2}{2}.
\]

A positive solution to \eqref{P} gives rise to a conformal metric $\widetilde g:=u^{4/(n-2)}g$ with respect to which $M$ has constant positive scalar curvature $c_n^{-1}a_n$ and its boundary $\partial M$ has constant mean curvature $d_n^{-1}b$, where $b$ is any given real number. The existence of a metric with these properties is one of the fundamental questions in conformal geometry on manifolds with boundary.

Our goal is to show the existence of a least-energy nodal solution for \eqref{P} when $b<0$. We begin with a brief overview, presenting the progress made to date, and we highlight the main ideas for addressing the problem at hand.
\smallskip

When $(M,g)$ is a closed manifold -- i.e., $\partial M=\emptyset$ -- the question boils down to proving the existence of a metric conformal to the given one for which the manifold has constant scalar curvature. This is the classical problem proposed and partially solved by Yamabe in 1960 \cite{Yamabe}; its definitive solution was not obtained until 1984, thanks to remarkable contributions by Trudinger, Aubin, and Schoen \cite{Trudinger,Aubin,Schoen}. A detailed account of this process can be found in the review article by Lee and Parker \cite{LeeParker}.

From a variational standpoint, the problem is equivalent to proving the existence of a positive, least-energy solution to the Yamabe equation
\begin{equation}\label{Yamabe}
-\Delta_g u+c_nR_gu=a_n|u|^{p-2}u\qquad\text{on }M
\end{equation}
or, equivalently, showing that the so-called \emph{Yamabe invariant} is attained. Aubin proved that if the Yamabe invariant for $(M,g)$ is strictly less than that of the standard sphere $\sn$, the Yamabe equation has a positive solution. For manifolds with a non-positive Yamabe invariant, this condition yields a solution immediately. When the Yamabe invariant is positive, the argument is more delicate. It relies on the construction of a suitable test function that depends on the manifold's geometry: for manifolds that are not locally conformally flat and have dimension $n\geq 6$, the leading correction is determined by the Weyl tensor; in the other cases, the argument involves the mass of the Green function and the positive mass theorem. Thus, the Weyl tensor and the mass are the two classical geometric quantities behind compactness.

The standard sphere is the model for the lack of compactness. Via stereographic projection, the least-energy positive solutions to the Yamabe equation on $\mathbb{S}^n$ correspond to the so-called \emph{standard bubbles}, which are the positive solutions to the Yamabe equation in $\mathbb{R}^n$. The invariance of this problem under dilations gives rise to a blow-up phenomenon responsible for the loss of compactness of \eqref{Yamabe}. This model remains fundamental in the problems described below.
\smallskip

For a compact manifold with boundary, the natural geometric question is whether there exists a metric, conformal to the given one, with prescribed constant scalar curvature in the interior and prescribed constant mean curvature on the boundary. As previously mentioned, in the case of positive scalar curvature, this is equivalent to proving the existence of a positive solution to \eqref{P}. In \cite{Cherrier}, Cherrier gave a condition for the existence of such a solution. Escobar proved its existence in a wide range of geometric contexts \cite{EscobarCPAM,EscobarAnnals,EscobarJDG,EscobarIndiana}, including two important limiting cases: the minimal boundary problem -- where the mean curvature of the boundary is prescribed to be zero -- and the scalar flat problem -- where the prescribed scalar curvature in the interior is zero. In both cases, the model domain is the standard hemisphere, and a condition for the existence of a positive solution is expressed via a strict inequality between the least possible energy of a solution for $M$ and the corresponding value for the hemisphere.

As with the problem for closed manifolds, the local geometry plays a role in the first non-trivial term of the energy expansion. At a nonumbilic boundary point, the trace-free second fundamental form provides the decisive correction; in other situations, the boundary Weyl tensor and a mass-type quantity enter the analysis. It is worth noting that the nonumbilicity condition is not artificial: for $n \ge 4$ it is robust and generic within natural classes of metrics, in the sense specified in \cite{CruzPistoia}. The complementary case -- in which all boundary points are umbilic -- is geometrically more rigid. This includes spherical caps, which are conformally equivalent to the Euclidean half-space.

For the problem where the prescribed interior scalar curvature is positive and the prescribed mean curvature of the boundary is an arbitrary constant -- that is, for problem \eqref{P} -- Han and Li \cite{HanLiDuke,HanLi} provided an existence criterion involving a particular spherical cap $B_b$ related to the given parameter $b$; see \eqref{eq:B_b}. It is expressed as
\begin{equation}\label{introHLthreshold}
\mu_1(M)<\mu_1(B_b),
\end{equation}
where $\mu_1(M)$ is the least possible energy of a non-trivial solution to \eqref{P} and $\mu_1(B_b)$ is the corresponding quantity for the spherical cap; see \eqref{eq:mu}.
For positive manifolds of dimension $n\geq 7$ with a nonumbilic boundary point, they proved that \eqref{introHLthreshold} holds for every $b\in\R$, thereby obtaining a positive least-energy solution to \eqref{P}. The restriction to the positive half-space of $\R^n$ of a particular bubble that depends on $b$ -- defined in \eqref{bubble} -- plays the role of the full bubble in the classical Yamabe problem and it also takes part in establishing the existence of a nodal solution. The complete classification of positive solutions in the Euclidean half-plane for any $b\in\r$ was given by Li and Zhu in \cite{LiZhu}.
\smallskip

Knowledge regarding sign-changing solutions remains limited. Leveraging its rich symmetry, Ding \cite{Ding} established the existence of infinitely many nodal solutions on the standard sphere. Other families of solutions have been constructed using gluing and symmetry methods; see, for example, \cite{delPinoEtAl2011,delPinoEtAl2013,MedinaMusso,Clapp2016,ClappFayaSaldaña}. Fernández and Petean \cite{FernandezPetean} obtained solutions on $\sn$ with a prescribed number of nodal domains.

For general closed manifolds, Ammann and Humbert \cite{AmmannHumbert} introduced a conformal invariant -- called the \emph{second Yamabe invariant} -- defined in terms of the second eigenvalue of the conformal Laplacian. They showed that if this invariant is attained by a function $u$, then $u$ changes sign, giving rise to a \emph{generalized metric} $|u|^{4/(n-2)}g$. In the aforementioned paper, they proved the existence of such a function $u$ on manifolds of dimension $n \ge 11$ that are not locally conformally flat. The existence of nodal solutions in lower dimensions has been established in specific cases, such as in \cite{ClappPistoiaTavares}. On the other hand, Premoselli and Vétois \cite{PremoselliVetois} proved that, on spheres of dimensions 3 through 10 with a metric close to the standard one, the second Yamabe invariant is not attained.
Thus, unlike the case of positive solutions, the existence of least-energy nodal solutions to the Yamabe problem on arbitrary closed manifolds is still far from being fully understood.
\smallskip

For manifolds with boundary, the available results on the existence of nodal solutions are quite recent. When the dimension is $n \ge 7$, Ho and Pyo \cite{HoPyo} addressed the minimal boundary problem -- that is, $b=0$ -- and established the existence of a nodal solution to \eqref{P} if the Weyl tensor is not identically zero on the boundary of $M$; while Ho, Lee, and Shin \cite{HoLeeShin} considered the scalar flat case -- obtained by setting $a_n=0$ and $b=1$ -- and obtained a nodal solution when the boundary is nonumbilic. Their strategy relied on extending the second Yamabe invariant to manifolds with boundary and showing that it is attained in the aforementioned cases.

The existence of a least-energy nodal solution for problem \eqref{P} with $b \ge 0$ was recently proved in \cite{CPP} for manifolds of dimension $n \ge 7$ whose boundary contains a nonumbilic point. In that work a criterion was established for the existence of a nodal solution with the minimum possible energy -- denoted by $\mu_2(M)$ -- stating that such a solution exists if
\begin{equation}\label{eq:mu2}
\mu_2(M) < \mu_1(M) + \mu_1(B_b),
\end{equation}
where $B_b$ is the spherical cap associated with the parameter $b$. The proof of this criterion relied on a minimax characterization of the level $\mu_2(M)$ and it utilized the fact that, for $b \ge 0$, there exists a neighborhood of the cones of positive and negative functions that is invariant under the negative gradient flow, as well as a result -- proven for $b \ge 0$ -- describing the loss of compactness below the level $\mu_1(M) + \mu_1(B_b)$. The test function introduced by Han and Li for manifolds with boundary, along with estimates obtained in \cite{HanLi}, were used to prove that \eqref{eq:mu2} is satisfied when the manifold has dimension $n \ge 7$ and its boundary contains a nonumbilic point. The restriction $b \ge 0$ plays a crucial role in the argument given to prove the invariance of the cones -- an argument that does not apply when $b < 0$, leaving this case open.

The spherical cap does not satisfy \eqref{eq:mu2}, as $\mu_2(B_b)=2 \mu_1(B_b)$ and, consequently, there is no least-energy nodal solution for \eqref{P} on $B_b$. Nevertheless, by exploiting its symmetries and employing a concentration-compactness argument, the existence of infinitely many nodal solutions was established in \cite{CPPcap} for $b\geq 0$ in dimensions $n=5$ and $n\ge7$, as well as in dimensions $n\geq 5$ when $b=0$. For values $b<0$ sufficiently close to $0$, multiple nodal solutions were also obtained in \cite{CPPcap}. Regarding the scalar flat problem with $b=1$, Almaraz and Wang \cite{AlmarazWang} found infinitely many nodal solutions on $B_b$ for dimensions $n\geq 4$, and gave precise information concerning their energy. These results, like those obtained for $\sn$, highlight the richness of nodal solutions at the thresholds of compactness.
\smallskip

This paper addresses the existence of a nodal solution to problem \eqref{P} for the remaining case $b<0$. The key observation is that invariance under the negative gradient flow of a neighborhood of the cones of positive and negative functions is not required; a localized deformation result near the Nehari nodal set suffices. We will present a more detailed exposition of the new strategy after stating our main result.

We assume throughout that $M$ is \emph{positive}, meaning that
\begin{equation}\label{positive}
\lambda_1(M):=\inf_{u\in H^1(M)\setminus\{0\}}
\frac{Q(u)}{\int_Mu^2\,dv_g}>0,
\end{equation}
where
\begin{equation}\label{Q}
Q(u):=\int_M\bigl(|\nabla_g u|^2+c_nR_gu^2\bigr)\,dv_g
+d_n\int_{\partial M}h_gu^2\,d\sigma_g.
\end{equation}
Our main result is the following.

\begin{theorem}\label{main}
Assume that $b<0$, $n\geq7$, and that $\partial M$ contains a nonumbilic point. Then, problem \eqref{P} has a least-energy sign-changing solution.
\end{theorem}

Combining Theorem \ref{main} with the result from \cite{CPP} provides a complete statement regarding the sign of the boundary parameter.

\begin{corollary}\label{allb}
Assume that $n\geq7$ and that $\partial M$ contains a nonumbilic point. Then, for every $b\in\R$, problem \eqref{P} has a least-energy sign-changing solution.
\end{corollary}

To the best of our knowledge, this is the first result establishing the existence of a least-energy nodal solution on a positive compact manifold with boundary, valid for all negative values  of the boundary parameter.

The proof consists of three steps. First, we use the minimax characterization of the 
level $\mu_2(M)$ given in \cite{CPP} and apply Willem's quantitative deformation 
lemma \cite[Lemma 2.3]{Willem} to generate a Palais--Smale sequence at the 
minimal nodal energy level $\mu_2(M)$, such that its distance from the cones of 
positive and negative functions exceeds a fixed positive constant; this result holds 
for all $b\in\mathbb{R}$ and does not require a neighborhood of the cones to be 
invariant under the negative gradient flow; see Proposition \ref{PSconstruction}. 
Second, we present a detailed analysis of the loss of compactness below the level 
$\mu_1(M)+\mu_1(B_b)$ for $b<0$. Based on these results, we obtain the same 
existence criterion as that established for $b\geq 0$ in \cite{CPP}; see Theorem 
\ref{abstract}. Finally, to verify the validity of this criterion under the hypotheses 
of Theorem \ref{main}, we use the test function of Han and Li, which is simpler for 
$b<0$ than for $b>0$. As they themselves noted, it suffices to consider the 
bubble associated with $b$ without the need for any correction term. The 
interactions between the positive minimizer and the bubble are of order 
$o(\eps^2)$ precisely when $n\ge7$.

We conclude this introduction by noting some interesting works, in the Euclidean context, that study the existence of a positive solution under more general critical boundary conditions. In \cite{PierottiTerracini95}, Pierroti and Terracini consider the case where the coefficient $b$ is not constant, and in \cite{PierottiTerracini97}, they study the situation where the linear part is replaced by a superlinear, but trace-subcritical, term.

The paper is organized as follows. In Section 2, we present the variational 
framework and collect the results from \cite{CPP} that are used here. In Section 3, 
we show the existence of a Palais--Smale sequence at the minimal nodal energy 
level that remains away from the cones. In Section 4, we describe the loss of 
compactness for $b<0$ at low energy levels, and in Section 5, we derive an 
abstract criterion for attaining the minimal nodal energy level. In Section 6 we 
prove the main theorem using the test function of Han and Li. Finally, in the 
Appendix, we prove an energy splitting that is used in Section 4.

\section{The variational setting}

Throughout this section, we assume that $b$ is any real number.

We set
$$\|u\|_Q:=\sqrt{Q(u)}$$
where $Q$ is defined in \eqref{Q}, and we take it as the norm in $H^1(M)$. Since $M$ is assumed to be positive, $\|\cdot\|_Q$ is equivalent to the standard norm, see \cite[Inequalities (2.2)]{CPP}. The inner product associated to it is
\begin{equation}
\langle u,v\rangle_Q:=\int_M\bigl(\nabla_gu\cdot\nabla_gv+c_nR_guv\bigr)\,dv_g
+d_n\int_{\partial M}h_guv\,d\sigma_g.
\end{equation}
The solutions to \eqref{P} are the critical points of the energy functional $J:H^1(M)\to\r$, given by
\begin{equation}\label{J}
J(u):=\frac12Q(u)-\frac{a_n}{p}\int_M|u|^p\,dv_g
-\frac{b}{q}\int_{\partial M}|u|^q\,d\sigma_g.
\end{equation}
It is of class $C^2$ and its derivative is
\begin{align}\label{Jprime}
J'(u)[v]=\langle u,v\rangle_Q-a_n\int_M|u|^{p-2}uv\,dv_g
-b\int_{\partial M}|u|^{q-2}uv\,d\sigma_g.
\end{align}
The nontrivial critical points of $J$ belong to the set
\[
\N:=\{u\in H^1(M)\setminus\{0\}:J'(u)[u]=0\}
\]
and those that change sign belong to the set
\[
\E:=\{u\in H^1(M):u^+\in\N,\ u^-\in\N\},
\]
where $u^+:=\max\{u,0\}$ and $u^-:=\min\{u,0\}$. We define
\begin{equation}\label{eq:mu}
\mu_1(M):=\inf_{u\in\N}J(u),\qquad
\mu_2(M):=\inf_{u\in\E}J(u).
\end{equation}
It is shown in \cite[Proposition 2.3]{CPP} that $\mu_1(M)$ and $\mu_2(M)$ are conformal invariants.

We shall use the following standard facts that were proven in \cite[Section 2]{CPP} for any $b\in\r$. The set $\N$ is a closed $C^2$ Hilbert submanifold of $H^1(M)$, called the Nehari manifold, and it is a natural constraint for $J$. For every $u\in H^1(M)$, $u\neq0$, there is a unique $t_u>0$ such that $t_uu\in\N$, and the function $t\mapsto J(tu)$ increases in $(0,t_u)$ and decreases in $(t_u,\infty)$. Therefore,
\begin{equation}\label{eq:t_u}
0=t_u\frac{d}{dt}J(tu)|_{t=t_u}=J'(t_uu)[t_uu],
\end{equation}
that is, $t_uu\in\N$. \ We will also use the fact that $\E$ is closed in $H^1(M)$. This is true because the maps $u\to u^+$ and $u\to u^-$ are continuous and $\N$ is closed in $H^1(M)$. We have $\E\subset\N$, and
$$\mu_2(M)>\mu_1(M)>0.$$
For later use we note that, if $b=-\beta<0$, then
\begin{align}
J(u)-\frac1pJ'(u)[u]
&=\left(\frac12-\frac1p\right)Q(u)
+\beta\left(\frac1q-\frac1p\right)\int_{\partial M}|u|^q\,d\sigma_g\notag\\
&=\frac1nQ(u)
+\beta\frac{n-2}{2n(n-1)}\int_{\partial M}|u|^q\,d\sigma_g.
\label{coerciveidentity}
\end{align}
In particular, if $J'(u_k)[u_k]=o(1)$ and $(u_k)$ is bounded away from zero, then its energy is asymptotically nonnegative.

We also recall the following minimax characterization of $\mu_2(M)$. Let
\[
\PP^+:=\{u\in H^1(M):u\geq0\},\qquad
\PP^-:=\{u\in H^1(M):u\leq0\}.
\]
Define $\tau:H^1(M)\to[0,\infty)$ to be $\tau(0):=0$ and, for $u\neq0$,
\begin{equation}\label{tau}
\tau(u):=
\begin{cases}
\frac{a_n\int_M|u|^p\,dv_g+b\int_{\partial M}|u|^q\,d\sigma_g}{Q(u)} &\text{if \ }b\geq 0, \\
\frac{a_n\int_M|u|^p\,dv_g}{Q(u)-b\int_{\partial M}|u|^q\,d\sigma_g}&\text{if \ }b< 0,
\end{cases}
\end{equation}
and let $\Sigma$ be the family of all continuous maps
$\sigma:[0,1]^2\to H^1(M)$ that satisfy
\begin{equation}\label{Sigma}
\begin{cases}
\sigma(0,t)\in\PP^+,&\sigma(1,t)\in\PP^-,\\
\sigma(s,0)=0,&J(\sigma(s,1))\leq0\text{ \ and \ }\tau(\sigma(s,1))\geq2.
\end{cases}
\end{equation}
Then, it is shown in \cite[Lemma 2.4]{CPP} that the image of every $\sigma\in\Sigma$ intersects $\E$ and that
\begin{equation}\label{minmax}
\mu_2(M)=\inf_{\sigma\in\Sigma}\,\max_{(s,t)\in[0,1]^2}J(\sigma(s,t)).
\end{equation}
Finally, as shown in \cite[Lemma 2.5$(a)$]{CPP}, there is $\alpha_0>0$ such that
\begin{equation}\label{separation}
B_{\alpha_0}(\PP^+)\cap\E=B_{\alpha_0}(\PP^-)\cap\E=\emptyset.
\end{equation}
where, for any $\mathcal{X}\subset H^1(M)$ and $\gamma>0$,
$$B_{\gamma}(\mathcal{X}):=\{u\in H^1(M):\dist_Q(u,\mathcal{X})\leq\gamma\}\qquad\text{and}\qquad \dist_Q(u,\mathcal{X}):=\inf_{v\in \mathcal{X}}\|u-v\|_Q.$$
It should be noted that the statements in \cite[Lemmas 2.4 and 2.5$(a)$]{CPP} do not impose any restriction on the sign of $b$: they hold for all $b\in\mathbb{R}$.

\section{A Palais--Smale sequence away from the positive and the negative cone}

The following proposition is the first new ingredient. For $b\geq 0$ it was proven in \cite[Lemma 2.6]{CPP} using that, in that case, $B_{\alpha_0}(\PP^+)$ and $B_{\alpha_0}(\PP^-)$ are strictly positively invariant under the negative gradient flow of $J$. The proof we give below does not use this invariance and applies to any $b \in\r$. It is based on Willem's quantitative deformation lemma \cite[Lemma 2.3]{Willem} and on the statements \eqref{minmax} and \eqref{separation}.

\begin{proposition}\label{PSconstruction}
There exists a sequence $(u_k)$ in $H^1(M)$ such that
\begin{equation}\label{PSnearE}
J(u_k)\to\mu_2(M),\qquad
\,  \|J'(u_k)\|_{H^{-1}(M)}\to0,\qquad
\dist_Q(u_k,\E)\to0.
\end{equation}
Consequently, for some $\alpha>0$ and all $k$ sufficiently large,
\begin{equation}\label{farcones}
\dist_Q(u_k,\PP^+\cup\PP^-)\geq\alpha.
\end{equation}
\end{proposition}

\begin{proof}
Let $\alpha_0$ be as in \eqref{separation} and assume there exist $\eps\in(0,\frac{\mu_2(M)}{2})$ and $\delta\in(0,\frac{\alpha_0}{3})$ such that
\begin{equation}\label{gradlower}
\|J'(u)\|_{H^{-1}(M)}\geq \frac{8\eps}{\delta}\qquad\text{for all \ }u\in J^{-1}[\mu_2(M)-2\eps,\mu_2(M)+2\eps]\cap B_{3\delta}(\E):=\mathcal{Z}.
\end{equation}
 We will show that this leads to a contradiction.

Indeed, by Willem's quantitative deformation lemma \cite[Lemma 2.3]{Willem}, there exists a continuous map $\eta:[0,1]\times H^1(M)\to H^1(M)$ such that
\begin{itemize}
\item[$(a)$] $\eta(r,u)=u$ if $r=0$ or $u\in H^1(M)\smallsetminus\mathcal{Z}$,
\item[$(b)$] $J(\eta(1,u))\leq \mu_2(M)-\eps$ if $u\in B_{\delta}(\Ecal)$ and $J(u)\leq \mu_2(M)+\eps$,
\item[$(c)$]
$\|\eta(t,u)-u\|_{Q}\leq \delta$, \, $\forall\,u\in H^{1}(M)$, $\forall \,t\in [0,1]$.
\end{itemize}
By the minimax characterization \eqref{minmax} of $\mu_2(M)$ there exists $\sigma\in\Sigma$ such that
\begin{equation}\label{eq:c+eps}
\max_{s,t\in[0,1]}J(\sigma(s,t))\leq \mu_2(M)+\eps.
\end{equation}
Define $\what\sigma:[0,1]\times[0,1]\to H^1(M)$ by
$$\what\sigma(s,t):=\eta(1,\sigma(s,t)).$$
Since $3\delta<\alpha_0$, statement \eqref{separation} implies that $\PP^+\cup\PP^-\subset H^1(M)\smallsetminus\mathcal{Z}$, so \eqref{Sigma} and item $(a)$ 
yields
$$\what\sigma(0,t)=\sigma(0,t)\in\PP^+,\qquad \what\sigma(1,t)=\sigma(1,t)\in\PP^-,\qquad \what\sigma(s,0)=\sigma(s,0)=0.$$
and, as $J(\sigma(s,1))\leq 0<\mu_2(M)-2\eps$, we get that $\what\sigma(s,1)=\sigma(s,1)$. Therefore,
$$J(\what\sigma(s,1))\leq 0\qquad\text{and}\qquad\tau(\what\sigma(s,1))\geq2.$$
This shows that $\what\sigma\in\Sigma$. 
On the other hand, by \cite[Lemma 2.4 (conclusion (b))]{CPP} there exists $(s_0,t_0)$ such that
$\what\sigma(s_0,t_0)\in\Ecal$, namely $\eta(1,\sigma(s_{0},t_{0})\in \Ecal$.
Set $u:=\sigma(s_0,t_0)$. Taking into account that $\eta(1,u)\in \Ecal$, and applying item $(c)$ we deduce that
\[
{\rm \dist}(u,\Ecal)\leq \|u-\eta(1,u)\|\leq \delta
\]
namely, $u\in B_\delta(\mathcal{E})$, and
$J(u)\le \mu_2(M)+\varepsilon$ by \eqref{eq:c+eps}. Hence, item $(b)$ yields
\[
J(\what\sigma(s_0,t_0))\le \mu_2(M)-\varepsilon<\mu_2(M)=\inf_{\mathcal{E}}J,
\]
contradicting the fact that $\what\sigma(s_0,t_0)\in\Ecal$.  

On the other hand, inequality \eqref{eq:c+eps} and item $(b)$ yield
$$
\max_{s,t\in[0,1]}J(\what\sigma(s,t))\leq \mu_2(M)-\eps.
$$
This contradicts the fact (shown in \cite[Lemma 2.4]{CPP}) that the image of any map in $\Sigma$ intersects $\E$, and proves that \eqref{gradlower} does not hold.

Consequently, by setting $\eps_k:=\frac{1}{k}$ and $\delta_k:=\sqrt{\eps_k}$, we conclude that there exists $u_k\in H^1(M)$ such that
$$J(u_k)\in[\mu_2(M)-2\eps_k,\mu_2(M)+2\eps_k],\qquad\|J'(u_k)\|_{H^-1(M)}<8\delta_k,\qquad \dist_Q(u_k,\E)\leq 3\delta_k,$$
for every $k$ sufficiently large. This proves the first statement.

To obtain \eqref{farcones}, choose $v_k\in\E$ with
$\|u_k-v_k\|_Q=o(1)$. By \eqref{separation},
\[
\dist_Q(u_k,\PP^\pm)
\geq\dist_Q(v_k,\PP^\pm)-\|u_k-v_k\|_Q
\geq\alpha_0-o(1).
\]
The statement is obtained by taking $\alpha:=\alpha_0/2$.
\end{proof}

\section{Compactness below the bubbling threshold for $b<0$}
\label{sec:compactness}

\begin{figure}[ht]
\centering
\begin{tikzpicture}[scale=0.82, line cap=round, line join=round]
\begin{scope}[xshift=0.0cm]
\shade[left color=planeL, right color=planeR, opacity=.9] (-1.107,-0.600) -- (1.188,-1.078) -- (3.401,0.122) -- (1.107,0.600) -- cycle;
\draw[planeR!70!black, thin] (-1.107,-0.600) -- (1.188,-1.078) -- (3.401,0.122) -- (1.107,0.600) -- cycle;
\draw[->, >=stealth] (0.000,0.000) -- (0.000,2.025) node[above] {$x_{n+1}$};
\draw[->, >=stealth] (0.000,0.000) -- (2.592,-0.540) node[below right=-1pt] {$x_n$};
\draw[->, >=stealth] (0.000,0.000) -- (-1.370,-0.743) node[below left=-2pt] {$x'$};
\node[blue!60!black] at (2.682,-0.045) {$\mathbb R^n_+$};
\draw[densely dotted, thin] (-1.164,-0.153) -- (-1.119,-0.162) -- (-1.073,-0.170) -- (-1.025,-0.177) -- (-0.976,-0.183) -- (-0.926,-0.189) -- (-0.876,-0.194) -- (-0.825,-0.197) -- (-0.773,-0.200) -- (-0.721,-0.202) -- (-0.669,-0.203) -- (-0.616,-0.203) -- (-0.564,-0.202) -- (-0.512,-0.201) -- (-0.460,-0.198) -- (-0.409,-0.194) -- (-0.358,-0.190) -- (-0.309,-0.185) -- (-0.260,-0.178) -- (-0.212,-0.171) -- (-0.165,-0.163) -- (-0.119,-0.155) -- (-0.075,-0.145) -- (-0.033,-0.135) -- (0.008,-0.124) -- (0.047,-0.112) -- (0.084,-0.100) -- (0.120,-0.087) -- (0.153,-0.073) -- (0.184,-0.059) -- (0.212,-0.044) -- (0.239,-0.029) -- (0.263,-0.013) -- (0.284,0.003) -- (0.303,0.019) -- (0.320,0.036) -- (0.334,0.053) -- (0.345,0.070) -- (0.353,0.087) -- (0.359,0.105) -- (0.362,0.122) -- (0.362,0.140) -- (0.360,0.157) -- (0.355,0.175) -- (0.347,0.192) -- (0.336,0.210) -- (0.323,0.227) -- (0.307,0.243) -- (0.288,0.260) -- (0.267,0.276) -- (0.244,0.292) -- (0.218,0.307) -- (0.190,0.322) -- (0.159,0.336) -- (0.126,0.350) -- (0.092,0.363) -- (0.055,0.376) -- (0.016,0.387) -- (-0.025,0.399) -- (-0.067,0.409) -- (-0.110,0.419) -- (-0.156,0.427) -- (-0.202,0.436) -- (-0.250,0.443) -- (-0.299,0.449) -- (-0.348,0.455) -- (-0.399,0.459) -- (-0.450,0.463) -- (-0.502,0.466) -- (-0.554,0.468) -- (-0.606,0.469) -- (-0.658,0.469) -- (-0.711,0.468) -- (-0.763,0.466) -- (-0.815,0.464) -- (-0.866,0.460) -- (-0.916,0.456) -- (-0.966,0.450) -- (-1.015,0.444) -- (-1.063,0.437) -- (-1.110,0.429) -- (-1.155,0.420) -- (-1.200,0.411) -- (-1.242,0.401) -- (-1.283,0.390) -- (-1.322,0.378) -- (-1.359,0.366) -- (-1.394,0.353) -- (-1.428,0.339) -- (-1.459,0.325) -- (-1.487,0.310) -- (-1.514,0.295) -- (-1.538,0.279) -- (-1.559,0.263) -- (-1.578,0.247) -- (-1.595,0.230) -- (-1.609,0.213) -- (-1.620,0.196) -- (-1.628,0.178) -- (-1.634,0.161) -- (-1.637,0.143) -- (-1.637,0.126) -- (-1.635,0.108) -- (-1.630,0.091) -- (-1.622,0.073) -- (-1.611,0.056) -- (-1.598,0.039) -- (-1.582,0.022) -- (-1.563,0.006) -- (-1.542,-0.010) -- (-1.519,-0.026) -- (-1.493,-0.041) -- (-1.464,-0.056) -- (-1.434,-0.070) -- (-1.401,-0.084) -- (-1.366,-0.097) -- (-1.330,-0.110) -- (-1.291,-0.122) -- (-1.250,-0.133) -- (-1.208,-0.143) -- (-1.164,-0.153) -- cycle;
\draw[dashed] (0.000,0.000) -- (-0.637,0.133) -- (-1.487,0.310);
\draw[capcol!60!black, densely dashed, thin] (-0.352,1.091) -- (-0.339,1.087) -- (-0.326,1.083) -- (-0.314,1.078) -- (-0.301,1.073) -- (-0.288,1.068) -- (-0.275,1.063) -- (-0.263,1.057) -- (-0.250,1.051) -- (-0.237,1.044) -- (-0.224,1.037) -- (-0.211,1.030) -- (-0.198,1.022) -- (-0.185,1.015) -- (-0.173,1.007) -- (-0.160,0.998) -- (-0.147,0.989) -- (-0.134,0.981) -- (-0.122,0.971) -- (-0.109,0.962) -- (-0.097,0.952) -- (-0.084,0.942) -- (-0.072,0.931) -- (-0.059,0.921) -- (-0.047,0.910) -- (-0.035,0.899) -- (-0.023,0.887) -- (-0.011,0.876) -- (0.001,0.864) -- (0.013,0.852) -- (0.025,0.840) -- (0.036,0.827) -- (0.048,0.814) -- (0.059,0.801) -- (0.070,0.788) -- (0.081,0.775) -- (0.092,0.761) -- (0.102,0.748) -- (0.113,0.734) -- (0.123,0.720) -- (0.133,0.706) -- (0.144,0.691) -- (0.153,0.677) -- (0.163,0.662) -- (0.172,0.647) -- (0.182,0.633) -- (0.191,0.618) -- (0.200,0.602) -- (0.208,0.587) -- (0.217,0.572) -- (0.225,0.557) -- (0.233,0.541) -- (0.241,0.525) -- (0.248,0.510) -- (0.256,0.494) -- (0.263,0.478) -- (0.270,0.463) -- (0.276,0.447) -- (0.283,0.431) -- (0.289,0.415) -- (0.295,0.399) -- (0.300,0.383) -- (0.306,0.367) -- (0.311,0.352) -- (0.316,0.336) -- (0.320,0.320) -- (0.325,0.304) -- (0.329,0.288) -- (0.333,0.273) -- (0.336,0.257) -- (0.340,0.241) -- (0.342,0.226) -- (0.345,0.210) -- (0.348,0.195) -- (0.350,0.180) -- (0.352,0.164) -- (0.353,0.149) -- (0.355,0.134) -- (0.356,0.120) -- (0.356,0.105) -- (0.357,0.090) -- (0.357,0.076) -- (0.357,0.062) -- (0.357,0.048) -- (0.356,0.034) -- (0.355,0.020) -- (0.354,0.006) -- (0.352,-0.007) -- (0.351,-0.020) -- (0.349,-0.034);
\draw[red!75!black, dashed] (-0.637,1.075) -- (0.000,0.000);
\shade[ball color=capcol, opacity=.85] (0.349,-0.034) -- (0.343,-0.062) -- (0.336,-0.089) -- (0.328,-0.115) -- (0.318,-0.140) -- (0.308,-0.163) -- (0.295,-0.186) -- (0.282,-0.207) -- (0.268,-0.226) -- (0.252,-0.244) -- (0.235,-0.261) -- (0.217,-0.276) -- (0.198,-0.289) -- (0.178,-0.301) -- (0.157,-0.312) -- (0.135,-0.320) -- (0.112,-0.327) -- (0.088,-0.332) -- (0.064,-0.336) -- (0.039,-0.338) -- (0.013,-0.338) -- (-0.013,-0.336) -- (-0.040,-0.333) -- (-0.067,-0.328) -- (-0.095,-0.321) -- (-0.123,-0.313) -- (-0.151,-0.303) -- (-0.179,-0.291) -- (-0.208,-0.278) -- (-0.236,-0.263) -- (-0.265,-0.246) -- (-0.293,-0.228) -- (-0.321,-0.209) -- (-0.350,-0.188) -- (-0.377,-0.166) -- (-0.405,-0.143) -- (-0.432,-0.118) -- (-0.458,-0.092) -- (-0.484,-0.065) -- (-0.510,-0.037) -- (-0.534,-0.008) -- (-0.559,0.022) -- (-0.582,0.053) -- (-0.604,0.085) -- (-0.626,0.117) -- (-0.646,0.150) -- (-0.666,0.184) -- (-0.684,0.218) -- (-0.702,0.252) -- (-0.718,0.287) -- (-0.733,0.322) -- (-0.747,0.358) -- (-0.760,0.393) -- (-0.772,0.428) -- (-0.782,0.463) -- (-0.791,0.498) -- (-0.799,0.533) -- (-0.805,0.567) -- (-0.810,0.601) -- (-0.813,0.635) -- (-0.815,0.668) -- (-0.816,0.700) -- (-0.815,0.731) -- (-0.813,0.762) -- (-0.810,0.791) -- (-0.805,0.820) -- (-0.799,0.848) -- (-0.792,0.875) -- (-0.783,0.900) -- (-0.773,0.924) -- (-0.761,0.947) -- (-0.748,0.969) -- (-0.735,0.989) -- (-0.719,1.008) -- (-0.703,1.026) -- (-0.686,1.042) -- (-0.667,1.056) -- (-0.648,1.069) -- (-0.627,1.080) -- (-0.606,1.089) -- (-0.584,1.097) -- (-0.560,1.103) -- (-0.536,1.108) -- (-0.512,1.111) -- (-0.486,1.112) -- (-0.460,1.111) -- (-0.434,1.109) -- (-0.407,1.104) -- (-0.379,1.099) -- (-0.352,1.091) -- (-0.352,1.091) -- (-0.340,1.088) -- (-0.328,1.084) -- (-0.317,1.080) -- (-0.305,1.076) -- (-0.294,1.072) -- (-0.283,1.068) -- (-0.271,1.063) -- (-0.260,1.059) -- (-0.249,1.054) -- (-0.237,1.049) -- (-0.226,1.044) -- (-0.215,1.039) -- (-0.204,1.034) -- (-0.193,1.029) -- (-0.182,1.023) -- (-0.172,1.018) -- (-0.161,1.012) -- (-0.150,1.006) -- (-0.140,1.000) -- (-0.129,0.994) -- (-0.119,0.988) -- (-0.108,0.981) -- (-0.098,0.975) -- (-0.088,0.968) -- (-0.078,0.962) -- (-0.068,0.955) -- (-0.058,0.948) -- (-0.048,0.941) -- (-0.038,0.933) -- (-0.028,0.926) -- (-0.019,0.919) -- (-0.009,0.911) -- (0.000,0.903) -- (0.010,0.895) -- (0.019,0.887) -- (0.028,0.879) -- (0.037,0.871) -- (0.046,0.863) -- (0.055,0.855) -- (0.063,0.846) -- (0.072,0.838) -- (0.081,0.829) -- (0.089,0.820) -- (0.097,0.811) -- (0.105,0.802) -- (0.114,0.793) -- (0.122,0.784) -- (0.129,0.775) -- (0.137,0.765) -- (0.145,0.756) -- (0.152,0.746) -- (0.160,0.737) -- (0.167,0.727) -- (0.174,0.717) -- (0.181,0.707) -- (0.188,0.697) -- (0.195,0.687) -- (0.202,0.677) -- (0.208,0.667) -- (0.215,0.656) -- (0.221,0.646) -- (0.227,0.635) -- (0.233,0.625) -- (0.239,0.614) -- (0.245,0.603) -- (0.251,0.593) -- (0.256,0.582) -- (0.262,0.571) -- (0.267,0.560) -- (0.272,0.549) -- (0.277,0.538) -- (0.282,0.527) -- (0.286,0.515) -- (0.291,0.504) -- (0.296,0.493) -- (0.300,0.481) -- (0.304,0.470) -- (0.308,0.459) -- (0.312,0.447) -- (0.316,0.435) -- (0.319,0.424) -- (0.323,0.412) -- (0.326,0.400) -- (0.329,0.389) -- (0.332,0.377) -- (0.335,0.365) -- (0.338,0.353) -- (0.341,0.341) -- (0.343,0.329) -- (0.345,0.318) -- (0.348,0.306) -- (0.350,0.294) -- (0.351,0.282) -- (0.353,0.269) -- (0.355,0.257) -- (0.356,0.245) -- (0.358,0.233) -- (0.359,0.221) -- (0.360,0.209) -- (0.361,0.197) -- (0.361,0.185) -- (0.362,0.173) -- (0.362,0.160) -- (0.362,0.148) -- (0.363,0.136) -- (0.363,0.124) -- (0.362,0.112) -- (0.362,0.100) -- (0.362,0.087) -- (0.361,0.075) -- (0.360,0.063) -- (0.359,0.051) -- (0.358,0.039) -- (0.357,0.027) -- (0.356,0.015) -- (0.354,0.003) -- (0.352,-0.009) -- (0.351,-0.022) -- (0.349,-0.034) -- cycle;
\draw[capcol!60!black, thick] (0.349,-0.034) -- (0.343,-0.062) -- (0.336,-0.089) -- (0.328,-0.115) -- (0.318,-0.140) -- (0.308,-0.163) -- (0.295,-0.186) -- (0.282,-0.207) -- (0.268,-0.226) -- (0.252,-0.244) -- (0.235,-0.261) -- (0.217,-0.276) -- (0.198,-0.289) -- (0.178,-0.301) -- (0.157,-0.312) -- (0.135,-0.320) -- (0.112,-0.327) -- (0.088,-0.332) -- (0.064,-0.336) -- (0.039,-0.338) -- (0.013,-0.338) -- (-0.013,-0.336) -- (-0.040,-0.333) -- (-0.067,-0.328) -- (-0.095,-0.321) -- (-0.123,-0.313) -- (-0.151,-0.303) -- (-0.179,-0.291) -- (-0.208,-0.278) -- (-0.236,-0.263) -- (-0.265,-0.246) -- (-0.293,-0.228) -- (-0.321,-0.209) -- (-0.350,-0.188) -- (-0.377,-0.166) -- (-0.405,-0.143) -- (-0.432,-0.118) -- (-0.458,-0.092) -- (-0.484,-0.065) -- (-0.510,-0.037) -- (-0.534,-0.008) -- (-0.559,0.022) -- (-0.582,0.053) -- (-0.604,0.085) -- (-0.626,0.117) -- (-0.646,0.150) -- (-0.666,0.184) -- (-0.684,0.218) -- (-0.702,0.252) -- (-0.718,0.287) -- (-0.733,0.322) -- (-0.747,0.358) -- (-0.760,0.393) -- (-0.772,0.428) -- (-0.782,0.463) -- (-0.791,0.498) -- (-0.799,0.533) -- (-0.805,0.567) -- (-0.810,0.601) -- (-0.813,0.635) -- (-0.815,0.668) -- (-0.816,0.700) -- (-0.815,0.731) -- (-0.813,0.762) -- (-0.810,0.791) -- (-0.805,0.820) -- (-0.799,0.848) -- (-0.792,0.875) -- (-0.783,0.900) -- (-0.773,0.924) -- (-0.761,0.947) -- (-0.748,0.969) -- (-0.735,0.989) -- (-0.719,1.008) -- (-0.703,1.026) -- (-0.686,1.042) -- (-0.667,1.056) -- (-0.648,1.069) -- (-0.627,1.080) -- (-0.606,1.089) -- (-0.584,1.097) -- (-0.560,1.103) -- (-0.536,1.108) -- (-0.512,1.111) -- (-0.486,1.112) -- (-0.460,1.111) -- (-0.434,1.109) -- (-0.407,1.104) -- (-0.379,1.099) -- (-0.352,1.091);
\draw (-0.637,0.133) circle (1);
\fill (0.000,0.000) circle (1.6pt);
\fill[red!75!black] (-0.637,1.075) circle (1.8pt);
\fill[red!75!black] (-0.637,0.133) circle (1.8pt);
\node[red!75!black, anchor=south east, font=\small] (etaL) at (-0.987,1.625) {$\eta_b$};
\draw[->, red!75!black, densely dotted, >=stealth, shorten >=2pt] (etaL.east) to[bend left=25] (-0.637,1.075);
\node[red!75!black, anchor=north east, font=\small] (xiL) at (-1.187,-1.017) {$\xi_b$};
\draw[->, red!75!black, densely dotted, >=stealth, shorten >=2pt] (xiL.north) to[bend left=30] (-0.637,0.133);
\node[red!75!black] at (0.353,1.188) {$B_b$};
\node[red!75!black] at (1.550,2.050) {$b>0$};
\end{scope}
\begin{scope}[xshift=6.0cm]
\shade[left color=planeL, right color=planeR, opacity=.9] (-1.107,-0.600) -- (1.188,-1.078) -- (3.401,0.122) -- (1.107,0.600) -- cycle;
\draw[planeR!70!black, thin] (-1.107,-0.600) -- (1.188,-1.078) -- (3.401,0.122) -- (1.107,0.600) -- cycle;
\draw[->, >=stealth] (0.000,0.000) -- (0.000,2.025) node[above] {$x_{n+1}$};
\draw[->, >=stealth] (0.000,0.000) -- (2.592,-0.540) node[below right=-1pt] {$x_n$};
\draw[->, >=stealth] (0.000,0.000) -- (-1.370,-0.743) node[below left=-2pt] {$x'$};
\node[blue!60!black] at (2.682,-0.045) {$\mathbb R^n_+$};
\draw[densely dotted, thin] (-0.527,-0.286) -- (-0.482,-0.295) -- (-0.435,-0.303) -- (-0.387,-0.310) -- (-0.339,-0.316) -- (-0.289,-0.322) -- (-0.239,-0.326) -- (-0.187,-0.330) -- (-0.136,-0.333) -- (-0.084,-0.335) -- (-0.031,-0.336) -- (0.021,-0.336) -- (0.073,-0.335) -- (0.125,-0.334) -- (0.177,-0.331) -- (0.228,-0.327) -- (0.279,-0.323) -- (0.329,-0.317) -- (0.378,-0.311) -- (0.426,-0.304) -- (0.473,-0.296) -- (0.518,-0.288) -- (0.562,-0.278) -- (0.605,-0.268) -- (0.645,-0.257) -- (0.685,-0.245) -- (0.722,-0.233) -- (0.757,-0.220) -- (0.790,-0.206) -- (0.821,-0.192) -- (0.850,-0.177) -- (0.876,-0.162) -- (0.900,-0.146) -- (0.922,-0.130) -- (0.941,-0.114) -- (0.957,-0.097) -- (0.971,-0.080) -- (0.982,-0.063) -- (0.991,-0.046) -- (0.996,-0.028) -- (1.000,-0.011) -- (1.000,0.007) -- (0.997,0.025) -- (0.992,0.042) -- (0.984,0.060) -- (0.974,0.077) -- (0.960,0.094) -- (0.944,0.111) -- (0.926,0.127) -- (0.905,0.143) -- (0.881,0.159) -- (0.855,0.174) -- (0.827,0.189) -- (0.797,0.203) -- (0.764,0.217) -- (0.729,0.230) -- (0.692,0.243) -- (0.653,0.254) -- (0.613,0.266) -- (0.571,0.276) -- (0.527,0.286) -- (0.482,0.295) -- (0.435,0.303) -- (0.387,0.310) -- (0.339,0.316) -- (0.289,0.322) -- (0.239,0.326) -- (0.187,0.330) -- (0.136,0.333) -- (0.084,0.335) -- (0.031,0.336) -- (-0.021,0.336) -- (-0.073,0.335) -- (-0.125,0.334) -- (-0.177,0.331) -- (-0.228,0.327) -- (-0.279,0.323) -- (-0.329,0.317) -- (-0.378,0.311) -- (-0.426,0.304) -- (-0.473,0.296) -- (-0.518,0.288) -- (-0.562,0.278) -- (-0.605,0.268) -- (-0.645,0.257) -- (-0.685,0.245) -- (-0.722,0.233) -- (-0.757,0.220) -- (-0.790,0.206) -- (-0.821,0.192) -- (-0.850,0.177) -- (-0.876,0.162) -- (-0.900,0.146) -- (-0.922,0.130) -- (-0.941,0.114) -- (-0.957,0.097) -- (-0.971,0.080) -- (-0.982,0.063) -- (-0.991,0.046) -- (-0.996,0.028) -- (-1.000,0.011) -- (-1.000,-0.007) -- (-0.997,-0.025) -- (-0.992,-0.042) -- (-0.984,-0.060) -- (-0.974,-0.077) -- (-0.960,-0.094) -- (-0.944,-0.111) -- (-0.926,-0.127) -- (-0.905,-0.143) -- (-0.881,-0.159) -- (-0.855,-0.174) -- (-0.827,-0.189) -- (-0.797,-0.203) -- (-0.764,-0.217) -- (-0.729,-0.230) -- (-0.692,-0.243) -- (-0.653,-0.254) -- (-0.613,-0.266) -- (-0.571,-0.276) -- (-0.527,-0.286) -- cycle;
\draw[capcol!60!black, densely dashed, thin] (0.204,0.979) -- (0.221,0.975) -- (0.238,0.969) -- (0.254,0.963) -- (0.270,0.955) -- (0.286,0.946) -- (0.302,0.936) -- (0.317,0.924) -- (0.331,0.912) -- (0.346,0.898) -- (0.359,0.884) -- (0.373,0.868) -- (0.386,0.851) -- (0.398,0.833) -- (0.410,0.814) -- (0.421,0.794) -- (0.432,0.773) -- (0.443,0.751) -- (0.453,0.728) -- (0.462,0.704) -- (0.470,0.679) -- (0.479,0.654) -- (0.486,0.627) -- (0.493,0.600) -- (0.499,0.572) -- (0.505,0.544) -- (0.510,0.514) -- (0.514,0.484) -- (0.518,0.454) -- (0.521,0.423) -- (0.524,0.391) -- (0.525,0.359) -- (0.526,0.327) -- (0.527,0.294) -- (0.527,0.260) -- (0.526,0.227) -- (0.524,0.193) -- (0.522,0.158) -- (0.520,0.124) -- (0.516,0.090) -- (0.512,0.055) -- (0.507,0.020) -- (0.502,-0.015) -- (0.496,-0.049) -- (0.489,-0.084) -- (0.482,-0.119) -- (0.474,-0.153) -- (0.466,-0.187) -- (0.457,-0.221) -- (0.447,-0.255) -- (0.437,-0.288) -- (0.427,-0.321) -- (0.416,-0.354) -- (0.404,-0.386) -- (0.392,-0.418) -- (0.379,-0.449) -- (0.366,-0.480) -- (0.352,-0.510) -- (0.338,-0.539) -- (0.324,-0.568) -- (0.309,-0.596) -- (0.293,-0.623) -- (0.278,-0.650) -- (0.262,-0.675) -- (0.246,-0.700) -- (0.229,-0.724) -- (0.212,-0.747) -- (0.195,-0.769) -- (0.178,-0.790) -- (0.160,-0.811) -- (0.142,-0.830) -- (0.124,-0.848) -- (0.106,-0.865) -- (0.088,-0.881) -- (0.069,-0.896) -- (0.051,-0.910) -- (0.032,-0.923) -- (0.014,-0.934) -- (-0.005,-0.944) -- (-0.024,-0.954) -- (-0.042,-0.962) -- (-0.061,-0.968) -- (-0.079,-0.974) -- (-0.097,-0.978) -- (-0.116,-0.982) -- (-0.134,-0.983) -- (-0.152,-0.984) -- (-0.169,-0.984) -- (-0.187,-0.982) -- (-0.204,-0.979);
\draw[red!75!black, dashed] (0.000,0.942) -- (0.000,0.000);
\shade[ball color=capcol, opacity=.85] (-0.204,-0.979) -- (-0.221,-0.975) -- (-0.238,-0.969) -- (-0.254,-0.963) -- (-0.270,-0.955) -- (-0.286,-0.946) -- (-0.302,-0.936) -- (-0.317,-0.924) -- (-0.331,-0.912) -- (-0.346,-0.898) -- (-0.359,-0.884) -- (-0.373,-0.868) -- (-0.386,-0.851) -- (-0.398,-0.833) -- (-0.410,-0.814) -- (-0.421,-0.794) -- (-0.432,-0.773) -- (-0.443,-0.751) -- (-0.453,-0.728) -- (-0.462,-0.704) -- (-0.470,-0.679) -- (-0.479,-0.654) -- (-0.486,-0.627) -- (-0.493,-0.600) -- (-0.499,-0.572) -- (-0.505,-0.544) -- (-0.510,-0.514) -- (-0.514,-0.484) -- (-0.518,-0.454) -- (-0.521,-0.423) -- (-0.524,-0.391) -- (-0.525,-0.359) -- (-0.526,-0.327) -- (-0.527,-0.294) -- (-0.527,-0.260) -- (-0.526,-0.227) -- (-0.524,-0.193) -- (-0.522,-0.158) -- (-0.520,-0.124) -- (-0.516,-0.090) -- (-0.512,-0.055) -- (-0.507,-0.020) -- (-0.502,0.015) -- (-0.496,0.049) -- (-0.489,0.084) -- (-0.482,0.119) -- (-0.474,0.153) -- (-0.466,0.187) -- (-0.457,0.221) -- (-0.447,0.255) -- (-0.437,0.288) -- (-0.427,0.321) -- (-0.416,0.354) -- (-0.404,0.386) -- (-0.392,0.418) -- (-0.379,0.449) -- (-0.366,0.480) -- (-0.352,0.510) -- (-0.338,0.539) -- (-0.324,0.568) -- (-0.309,0.596) -- (-0.293,0.623) -- (-0.278,0.650) -- (-0.262,0.675) -- (-0.246,0.700) -- (-0.229,0.724) -- (-0.212,0.747) -- (-0.195,0.769) -- (-0.178,0.790) -- (-0.160,0.811) -- (-0.142,0.830) -- (-0.124,0.848) -- (-0.106,0.865) -- (-0.088,0.881) -- (-0.069,0.896) -- (-0.051,0.910) -- (-0.032,0.923) -- (-0.014,0.934) -- (0.005,0.944) -- (0.024,0.954) -- (0.042,0.962) -- (0.061,0.968) -- (0.079,0.974) -- (0.097,0.978) -- (0.116,0.982) -- (0.134,0.983) -- (0.152,0.984) -- (0.169,0.984) -- (0.187,0.982) -- (0.204,0.979) -- (0.204,0.979) -- (0.230,0.973) -- (0.255,0.967) -- (0.281,0.960) -- (0.306,0.952) -- (0.331,0.944) -- (0.356,0.935) -- (0.380,0.925) -- (0.405,0.914) -- (0.429,0.903) -- (0.452,0.892) -- (0.476,0.880) -- (0.499,0.867) -- (0.522,0.853) -- (0.544,0.839) -- (0.566,0.824) -- (0.587,0.809) -- (0.609,0.793) -- (0.629,0.777) -- (0.650,0.760) -- (0.669,0.743) -- (0.689,0.725) -- (0.708,0.706) -- (0.726,0.688) -- (0.744,0.668) -- (0.761,0.648) -- (0.778,0.628) -- (0.795,0.607) -- (0.810,0.586) -- (0.826,0.564) -- (0.840,0.542) -- (0.854,0.520) -- (0.868,0.497) -- (0.880,0.474) -- (0.893,0.451) -- (0.904,0.427) -- (0.915,0.403) -- (0.925,0.379) -- (0.935,0.354) -- (0.944,0.329) -- (0.953,0.304) -- (0.960,0.279) -- (0.967,0.254) -- (0.974,0.228) -- (0.979,0.202) -- (0.984,0.176) -- (0.989,0.150) -- (0.992,0.124) -- (0.995,0.098) -- (0.997,0.072) -- (0.999,0.045) -- (1.000,0.019) -- (1.000,-0.007) -- (0.999,-0.034) -- (0.998,-0.060) -- (0.996,-0.087) -- (0.994,-0.113) -- (0.990,-0.139) -- (0.986,-0.165) -- (0.982,-0.191) -- (0.976,-0.217) -- (0.970,-0.243) -- (0.963,-0.268) -- (0.956,-0.294) -- (0.948,-0.319) -- (0.939,-0.344) -- (0.930,-0.368) -- (0.920,-0.393) -- (0.909,-0.417) -- (0.898,-0.441) -- (0.886,-0.464) -- (0.873,-0.487) -- (0.860,-0.510) -- (0.846,-0.533) -- (0.832,-0.555) -- (0.817,-0.577) -- (0.801,-0.598) -- (0.785,-0.619) -- (0.769,-0.640) -- (0.752,-0.660) -- (0.734,-0.679) -- (0.716,-0.698) -- (0.697,-0.717) -- (0.678,-0.735) -- (0.658,-0.753) -- (0.638,-0.770) -- (0.618,-0.786) -- (0.597,-0.803) -- (0.575,-0.818) -- (0.553,-0.833) -- (0.531,-0.847) -- (0.509,-0.861) -- (0.486,-0.874) -- (0.463,-0.887) -- (0.439,-0.898) -- (0.415,-0.910) -- (0.391,-0.920) -- (0.367,-0.930) -- (0.342,-0.940) -- (0.317,-0.948) -- (0.292,-0.956) -- (0.266,-0.964) -- (0.241,-0.971) -- (0.215,-0.977) -- (0.189,-0.982) -- (0.163,-0.987) -- (0.137,-0.991) -- (0.111,-0.994) -- (0.085,-0.996) -- (0.058,-0.998) -- (0.032,-0.999) -- (0.006,-1.000) -- (-0.021,-1.000) -- (-0.047,-0.999) -- (-0.073,-0.997) -- (-0.100,-0.995) -- (-0.126,-0.992) -- (-0.152,-0.988) -- (-0.178,-0.984) -- (-0.204,-0.979) -- cycle;
\draw[capcol!60!black, thick] (-0.204,-0.979) -- (-0.221,-0.975) -- (-0.238,-0.969) -- (-0.254,-0.963) -- (-0.270,-0.955) -- (-0.286,-0.946) -- (-0.302,-0.936) -- (-0.317,-0.924) -- (-0.331,-0.912) -- (-0.346,-0.898) -- (-0.359,-0.884) -- (-0.373,-0.868) -- (-0.386,-0.851) -- (-0.398,-0.833) -- (-0.410,-0.814) -- (-0.421,-0.794) -- (-0.432,-0.773) -- (-0.443,-0.751) -- (-0.453,-0.728) -- (-0.462,-0.704) -- (-0.470,-0.679) -- (-0.479,-0.654) -- (-0.486,-0.627) -- (-0.493,-0.600) -- (-0.499,-0.572) -- (-0.505,-0.544) -- (-0.510,-0.514) -- (-0.514,-0.484) -- (-0.518,-0.454) -- (-0.521,-0.423) -- (-0.524,-0.391) -- (-0.525,-0.359) -- (-0.526,-0.327) -- (-0.527,-0.294) -- (-0.527,-0.260) -- (-0.526,-0.227) -- (-0.524,-0.193) -- (-0.522,-0.158) -- (-0.520,-0.124) -- (-0.516,-0.090) -- (-0.512,-0.055) -- (-0.507,-0.020) -- (-0.502,0.015) -- (-0.496,0.049) -- (-0.489,0.084) -- (-0.482,0.119) -- (-0.474,0.153) -- (-0.466,0.187) -- (-0.457,0.221) -- (-0.447,0.255) -- (-0.437,0.288) -- (-0.427,0.321) -- (-0.416,0.354) -- (-0.404,0.386) -- (-0.392,0.418) -- (-0.379,0.449) -- (-0.366,0.480) -- (-0.352,0.510) -- (-0.338,0.539) -- (-0.324,0.568) -- (-0.309,0.596) -- (-0.293,0.623) -- (-0.278,0.650) -- (-0.262,0.675) -- (-0.246,0.700) -- (-0.229,0.724) -- (-0.212,0.747) -- (-0.195,0.769) -- (-0.178,0.790) -- (-0.160,0.811) -- (-0.142,0.830) -- (-0.124,0.848) -- (-0.106,0.865) -- (-0.088,0.881) -- (-0.069,0.896) -- (-0.051,0.910) -- (-0.032,0.923) -- (-0.014,0.934) -- (0.005,0.944) -- (0.024,0.954) -- (0.042,0.962) -- (0.061,0.968) -- (0.079,0.974) -- (0.097,0.978) -- (0.116,0.982) -- (0.134,0.983) -- (0.152,0.984) -- (0.169,0.984) -- (0.187,0.982) -- (0.204,0.979);
\draw (0.000,0.000) circle (1);
\fill (0.000,0.000) circle (1.6pt);
\fill[red!75!black] (0.000,0.942) circle (1.8pt);
\fill[red!75!black] (0.000,0.000) circle (1.8pt);
\node[red!75!black, anchor=south east, font=\small] (etaL) at (-0.450,1.342) {$\eta_b$};
\draw[->, red!75!black, densely dotted, >=stealth, shorten >=2pt] (etaL.east) to[bend left=25] (0.000,0.942);
\node[red!75!black, anchor=north east, font=\small] (xiL) at (-0.750,-1.250) {$\xi_b$};
\draw[->, red!75!black, densely dotted, >=stealth, shorten >=2pt] (xiL.north) to[bend left=30] (0.000,0.000);
\node[red!75!black] (BL) at (0.900,1.350) {$B_b$};
\draw[->, red!75!black, >=stealth] (BL.south west) -- (0.350,0.600);
\node[red!75!black] at (1.550,2.050) {$b=0$};
\end{scope}
\begin{scope}[xshift=12.0cm]
\shade[left color=planeL, right color=planeR, opacity=.9] (-1.107,-0.600) -- (1.188,-1.078) -- (3.401,0.122) -- (1.107,0.600) -- cycle;
\draw[planeR!70!black, thin] (-1.107,-0.600) -- (1.188,-1.078) -- (3.401,0.122) -- (1.107,0.600) -- cycle;
\draw[->, >=stealth] (0.000,0.000) -- (0.000,2.025) node[above] {$x_{n+1}$};
\draw[->, >=stealth] (0.000,0.000) -- (2.592,-0.540) node[below right=-1pt] {$x_n$};
\draw[->, >=stealth] (0.000,0.000) -- (-1.370,-0.743) node[below left=-2pt] {$x'$};
\node[blue!60!black] at (2.682,-0.045) {$\mathbb R^n_+$};
\draw[densely dotted, thin] (0.110,-0.419) -- (0.156,-0.427) -- (0.202,-0.436) -- (0.250,-0.443) -- (0.299,-0.449) -- (0.348,-0.455) -- (0.399,-0.459) -- (0.450,-0.463) -- (0.502,-0.466) -- (0.554,-0.468) -- (0.606,-0.469) -- (0.658,-0.469) -- (0.711,-0.468) -- (0.763,-0.466) -- (0.815,-0.464) -- (0.866,-0.460) -- (0.916,-0.456) -- (0.966,-0.450) -- (1.015,-0.444) -- (1.063,-0.437) -- (1.110,-0.429) -- (1.155,-0.420) -- (1.200,-0.411) -- (1.242,-0.401) -- (1.283,-0.390) -- (1.322,-0.378) -- (1.359,-0.366) -- (1.394,-0.353) -- (1.428,-0.339) -- (1.459,-0.325) -- (1.487,-0.310) -- (1.514,-0.295) -- (1.538,-0.279) -- (1.559,-0.263) -- (1.578,-0.247) -- (1.595,-0.230) -- (1.609,-0.213) -- (1.620,-0.196) -- (1.628,-0.178) -- (1.634,-0.161) -- (1.637,-0.143) -- (1.637,-0.126) -- (1.635,-0.108) -- (1.630,-0.091) -- (1.622,-0.073) -- (1.611,-0.056) -- (1.598,-0.039) -- (1.582,-0.022) -- (1.563,-0.006) -- (1.542,0.010) -- (1.519,0.026) -- (1.493,0.041) -- (1.464,0.056) -- (1.434,0.070) -- (1.401,0.084) -- (1.366,0.097) -- (1.330,0.110) -- (1.291,0.122) -- (1.250,0.133) -- (1.208,0.143) -- (1.164,0.153) -- (1.119,0.162) -- (1.073,0.170) -- (1.025,0.177) -- (0.976,0.183) -- (0.926,0.189) -- (0.876,0.194) -- (0.825,0.197) -- (0.773,0.200) -- (0.721,0.202) -- (0.669,0.203) -- (0.616,0.203) -- (0.564,0.202) -- (0.512,0.201) -- (0.460,0.198) -- (0.409,0.194) -- (0.358,0.190) -- (0.309,0.185) -- (0.260,0.178) -- (0.212,0.171) -- (0.165,0.163) -- (0.119,0.155) -- (0.075,0.145) -- (0.033,0.135) -- (-0.008,0.124) -- (-0.047,0.112) -- (-0.084,0.100) -- (-0.120,0.087) -- (-0.153,0.073) -- (-0.184,0.059) -- (-0.212,0.044) -- (-0.239,0.029) -- (-0.263,0.013) -- (-0.284,-0.003) -- (-0.303,-0.019) -- (-0.320,-0.036) -- (-0.334,-0.053) -- (-0.345,-0.070) -- (-0.353,-0.087) -- (-0.359,-0.105) -- (-0.362,-0.122) -- (-0.362,-0.140) -- (-0.360,-0.157) -- (-0.355,-0.175) -- (-0.347,-0.192) -- (-0.336,-0.210) -- (-0.323,-0.227) -- (-0.307,-0.243) -- (-0.288,-0.260) -- (-0.267,-0.276) -- (-0.244,-0.292) -- (-0.218,-0.307) -- (-0.190,-0.322) -- (-0.159,-0.336) -- (-0.126,-0.350) -- (-0.092,-0.363) -- (-0.055,-0.376) -- (-0.016,-0.387) -- (0.025,-0.399) -- (0.067,-0.409) -- (0.110,-0.419) -- cycle;
\draw[capcol!60!black, densely dashed, thin] (0.783,0.856) -- (0.790,0.855) -- (0.797,0.853) -- (0.802,0.849) -- (0.807,0.845) -- (0.810,0.840) -- (0.813,0.835) -- (0.815,0.828) -- (0.816,0.821) -- (0.816,0.813) -- (0.815,0.804) -- (0.814,0.794) -- (0.811,0.784) -- (0.808,0.773) -- (0.804,0.761) -- (0.799,0.748) -- (0.793,0.735) -- (0.786,0.721) -- (0.778,0.707) -- (0.770,0.692) -- (0.761,0.676) -- (0.751,0.659) -- (0.740,0.643) -- (0.728,0.625) -- (0.716,0.607) -- (0.703,0.589) -- (0.689,0.570) -- (0.675,0.551) -- (0.660,0.531) -- (0.644,0.511) -- (0.628,0.491) -- (0.611,0.470) -- (0.593,0.449) -- (0.575,0.428) -- (0.557,0.407) -- (0.538,0.385) -- (0.518,0.364) -- (0.498,0.342) -- (0.478,0.320) -- (0.457,0.298) -- (0.436,0.276) -- (0.415,0.254) -- (0.393,0.233) -- (0.371,0.211) -- (0.349,0.189) -- (0.327,0.168) -- (0.304,0.147) -- (0.282,0.126) -- (0.259,0.105) -- (0.236,0.084) -- (0.214,0.064) -- (0.191,0.044) -- (0.168,0.025) -- (0.146,0.006) -- (0.124,-0.013) -- (0.101,-0.031) -- (0.079,-0.049) -- (0.057,-0.066) -- (0.036,-0.083) -- (0.015,-0.099) -- (-0.006,-0.115) -- (-0.027,-0.130) -- (-0.047,-0.144) -- (-0.067,-0.158) -- (-0.086,-0.171) -- (-0.105,-0.183) -- (-0.124,-0.194) -- (-0.141,-0.205) -- (-0.159,-0.215) -- (-0.175,-0.225) -- (-0.191,-0.233) -- (-0.207,-0.241) -- (-0.222,-0.248) -- (-0.236,-0.254) -- (-0.249,-0.260) -- (-0.262,-0.264) -- (-0.274,-0.268) -- (-0.285,-0.271) -- (-0.296,-0.273) -- (-0.306,-0.274) -- (-0.314,-0.275) -- (-0.323,-0.274) -- (-0.330,-0.273) -- (-0.336,-0.271) -- (-0.342,-0.268) -- (-0.347,-0.264) -- (-0.350,-0.259) -- (-0.353,-0.254) -- (-0.356,-0.247) -- (-0.357,-0.240);
\draw[dashed] (0.000,0.000) -- (0.637,-0.133);
\draw[red!75!black, dashed] (0.637,0.809) -- (0.000,0.000);
\shade[ball color=capcol, opacity=.85] (-0.357,-0.240) -- (-0.357,-0.234) -- (-0.357,-0.227) -- (-0.356,-0.219) -- (-0.355,-0.211) -- (-0.353,-0.203) -- (-0.350,-0.194) -- (-0.347,-0.184) -- (-0.343,-0.174) -- (-0.339,-0.164) -- (-0.334,-0.153) -- (-0.328,-0.141) -- (-0.322,-0.129) -- (-0.315,-0.117) -- (-0.308,-0.104) -- (-0.300,-0.091) -- (-0.292,-0.078) -- (-0.283,-0.064) -- (-0.274,-0.049) -- (-0.264,-0.035) -- (-0.254,-0.020) -- (-0.243,-0.005) -- (-0.232,0.011) -- (-0.220,0.027) -- (-0.207,0.043) -- (-0.195,0.059) -- (-0.182,0.076) -- (-0.168,0.093) -- (-0.154,0.110) -- (-0.140,0.127) -- (-0.125,0.144) -- (-0.110,0.162) -- (-0.094,0.179) -- (-0.079,0.197) -- (-0.063,0.215) -- (-0.046,0.233) -- (-0.029,0.251) -- (-0.013,0.269) -- (0.005,0.287) -- (0.022,0.305) -- (0.040,0.323) -- (0.058,0.341) -- (0.076,0.359) -- (0.094,0.377) -- (0.112,0.395) -- (0.130,0.413) -- (0.149,0.430) -- (0.167,0.448) -- (0.186,0.465) -- (0.205,0.482) -- (0.224,0.499) -- (0.242,0.516) -- (0.261,0.532) -- (0.280,0.548) -- (0.298,0.564) -- (0.317,0.580) -- (0.335,0.595) -- (0.354,0.610) -- (0.372,0.625) -- (0.390,0.640) -- (0.408,0.654) -- (0.426,0.667) -- (0.443,0.681) -- (0.461,0.694) -- (0.478,0.706) -- (0.495,0.718) -- (0.511,0.730) -- (0.527,0.741) -- (0.543,0.752) -- (0.559,0.762) -- (0.574,0.772) -- (0.589,0.781) -- (0.604,0.790) -- (0.618,0.798) -- (0.632,0.806) -- (0.645,0.813) -- (0.658,0.820) -- (0.671,0.826) -- (0.683,0.832) -- (0.695,0.837) -- (0.706,0.841) -- (0.717,0.845) -- (0.727,0.849) -- (0.736,0.852) -- (0.746,0.854) -- (0.754,0.855) -- (0.762,0.857) -- (0.770,0.857) -- (0.777,0.857) -- (0.783,0.856) -- (0.783,0.856) -- (0.820,0.850) -- (0.857,0.843) -- (0.893,0.834) -- (0.929,0.824) -- (0.965,0.812) -- (1.000,0.799) -- (1.035,0.785) -- (1.069,0.769) -- (1.102,0.752) -- (1.135,0.734) -- (1.167,0.715) -- (1.199,0.695) -- (1.229,0.673) -- (1.259,0.650) -- (1.288,0.627) -- (1.316,0.602) -- (1.343,0.576) -- (1.369,0.549) -- (1.394,0.521) -- (1.418,0.492) -- (1.441,0.462) -- (1.463,0.432) -- (1.483,0.401) -- (1.503,0.369) -- (1.521,0.336) -- (1.538,0.302) -- (1.553,0.268) -- (1.568,0.234) -- (1.581,0.199) -- (1.593,0.163) -- (1.603,0.127) -- (1.612,0.091) -- (1.620,0.054) -- (1.626,0.017) -- (1.631,-0.020) -- (1.635,-0.057) -- (1.637,-0.095) -- (1.637,-0.132) -- (1.637,-0.170) -- (1.635,-0.207) -- (1.631,-0.244) -- (1.626,-0.281) -- (1.620,-0.318) -- (1.612,-0.355) -- (1.603,-0.391) -- (1.593,-0.427) -- (1.581,-0.463) -- (1.568,-0.498) -- (1.554,-0.533) -- (1.538,-0.567) -- (1.521,-0.600) -- (1.503,-0.633) -- (1.484,-0.665) -- (1.463,-0.696) -- (1.442,-0.727) -- (1.419,-0.757) -- (1.395,-0.786) -- (1.370,-0.813) -- (1.344,-0.840) -- (1.317,-0.866) -- (1.289,-0.891) -- (1.260,-0.915) -- (1.231,-0.938) -- (1.200,-0.960) -- (1.169,-0.980) -- (1.136,-0.999) -- (1.104,-1.018) -- (1.070,-1.034) -- (1.036,-1.050) -- (1.002,-1.064) -- (0.966,-1.077) -- (0.931,-1.089) -- (0.895,-1.099) -- (0.858,-1.108) -- (0.822,-1.116) -- (0.785,-1.122) -- (0.748,-1.127) -- (0.710,-1.130) -- (0.673,-1.132) -- (0.635,-1.133) -- (0.598,-1.132) -- (0.561,-1.130) -- (0.523,-1.126) -- (0.486,-1.121) -- (0.449,-1.115) -- (0.413,-1.107) -- (0.376,-1.098) -- (0.340,-1.088) -- (0.305,-1.076) -- (0.270,-1.063) -- (0.235,-1.048) -- (0.201,-1.033) -- (0.168,-1.016) -- (0.135,-0.997) -- (0.103,-0.978) -- (0.072,-0.957) -- (0.041,-0.936) -- (0.011,-0.913) -- (-0.017,-0.889) -- (-0.045,-0.864) -- (-0.072,-0.838) -- (-0.098,-0.811) -- (-0.123,-0.783) -- (-0.147,-0.754) -- (-0.169,-0.724) -- (-0.191,-0.693) -- (-0.211,-0.662) -- (-0.230,-0.630) -- (-0.248,-0.597) -- (-0.265,-0.563) -- (-0.281,-0.529) -- (-0.295,-0.494) -- (-0.308,-0.459) -- (-0.319,-0.424) -- (-0.330,-0.388) -- (-0.338,-0.351) -- (-0.346,-0.314) -- (-0.352,-0.278) -- (-0.357,-0.240) -- cycle;
\draw[capcol!60!black, thick] (-0.357,-0.240) -- (-0.357,-0.234) -- (-0.357,-0.227) -- (-0.356,-0.219) -- (-0.355,-0.211) -- (-0.353,-0.203) -- (-0.350,-0.194) -- (-0.347,-0.184) -- (-0.343,-0.174) -- (-0.339,-0.164) -- (-0.334,-0.153) -- (-0.328,-0.141) -- (-0.322,-0.129) -- (-0.315,-0.117) -- (-0.308,-0.104) -- (-0.300,-0.091) -- (-0.292,-0.078) -- (-0.283,-0.064) -- (-0.274,-0.049) -- (-0.264,-0.035) -- (-0.254,-0.020) -- (-0.243,-0.005) -- (-0.232,0.011) -- (-0.220,0.027) -- (-0.207,0.043) -- (-0.195,0.059) -- (-0.182,0.076) -- (-0.168,0.093) -- (-0.154,0.110) -- (-0.140,0.127) -- (-0.125,0.144) -- (-0.110,0.162) -- (-0.094,0.179) -- (-0.079,0.197) -- (-0.063,0.215) -- (-0.046,0.233) -- (-0.029,0.251) -- (-0.013,0.269) -- (0.005,0.287) -- (0.022,0.305) -- (0.040,0.323) -- (0.058,0.341) -- (0.076,0.359) -- (0.094,0.377) -- (0.112,0.395) -- (0.130,0.413) -- (0.149,0.430) -- (0.167,0.448) -- (0.186,0.465) -- (0.205,0.482) -- (0.224,0.499) -- (0.242,0.516) -- (0.261,0.532) -- (0.280,0.548) -- (0.298,0.564) -- (0.317,0.580) -- (0.335,0.595) -- (0.354,0.610) -- (0.372,0.625) -- (0.390,0.640) -- (0.408,0.654) -- (0.426,0.667) -- (0.443,0.681) -- (0.461,0.694) -- (0.478,0.706) -- (0.495,0.718) -- (0.511,0.730) -- (0.527,0.741) -- (0.543,0.752) -- (0.559,0.762) -- (0.574,0.772) -- (0.589,0.781) -- (0.604,0.790) -- (0.618,0.798) -- (0.632,0.806) -- (0.645,0.813) -- (0.658,0.820) -- (0.671,0.826) -- (0.683,0.832) -- (0.695,0.837) -- (0.706,0.841) -- (0.717,0.845) -- (0.727,0.849) -- (0.736,0.852) -- (0.746,0.854) -- (0.754,0.855) -- (0.762,0.857) -- (0.770,0.857) -- (0.777,0.857) -- (0.783,0.856);
\draw (0.637,-0.133) circle (1);
\fill (0.000,0.000) circle (1.6pt);
\fill[red!75!black] (0.637,0.809) circle (1.8pt);
\fill[red!75!black] (0.637,-0.133) circle (1.8pt);
\node[red!75!black, anchor=south east, font=\small] (etaL) at (0.187,1.209) {$\eta_b$};
\draw[->, red!75!black, densely dotted, >=stealth, shorten >=2pt] (etaL.east) to[bend left=25] (0.637,0.809);
\node[red!75!black, anchor=north east, font=\small] (xiL) at (-0.113,-1.383) {$\xi_b$};
\draw[->, red!75!black, densely dotted, >=stealth, shorten >=2pt] (xiL.north) to[bend left=30] (0.637,-0.133);
\node[red!75!black] (BL) at (1.537,1.217) {$B_b$};
\draw[->, red!75!black, >=stealth] (BL.south west) -- (0.987,0.467);
\node[red!75!black] at (1.550,2.050) {$b<0$};
\end{scope}
\end{tikzpicture}
\caption{\small{The spherical cap $B_b=\pi_b^{-1}(\rn_+)$ for $b>0$ (left), $b=0$ (center) and $b<0$ (right). Here $\xi_b=(0,\ldots,0,-\frac{b}{n-2},0)$ is the center of the unit sphere, $\eta_b=(0,\ldots,0,-\frac{b}{n-2},1)$ is the north pole and $x'=(x_1,\ldots,x_{n-1})$.}}
\label{fig:caps}
\end{figure}
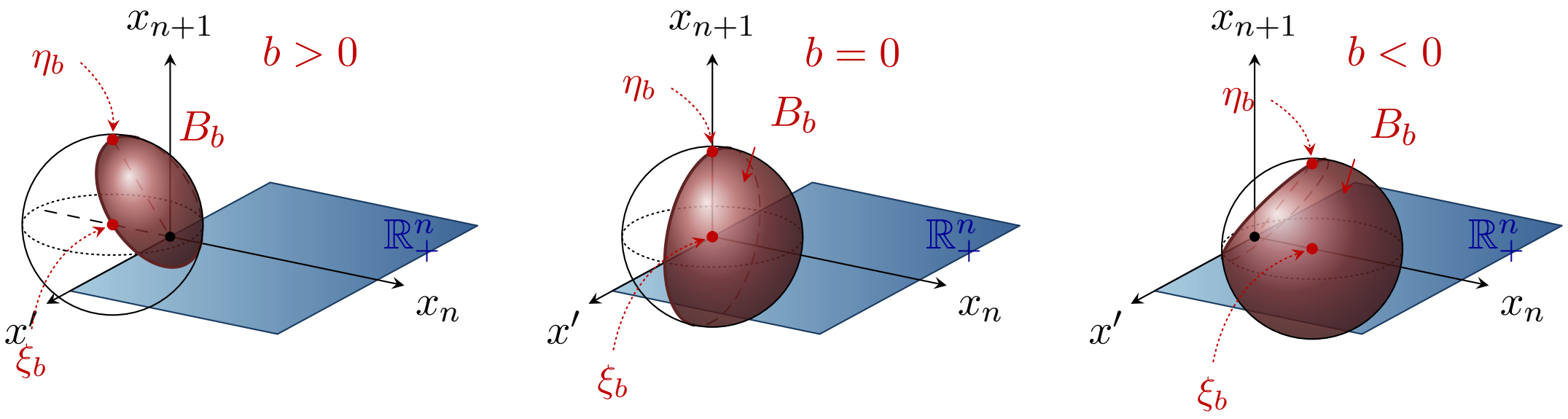

For any $b\in\r$, let $B_b$ be the spherical cap
\begin{equation}\label{eq:B_b}
B_b:=\pi_b^{-1}(\rn_+),
\end{equation}
where $\rn_+:=\{(x_1,\ldots,x_n)\in\rn:x_n\geq 0\}$ and $\pi_b$ is the stereographic projection from the north pole $\eta_b:=(0,\ldots,0,-\frac{b}{n-2},1)$ of the sphere of radius $1$ in $\r^{n+1}$ centered at $\xi_b:=(0,\ldots,0,-\frac{b}{n-2},0)$; see Figure~\ref{fig:caps}. Given that the center $\xi_b$ lies on the positive $x_n$-axis when $b<0$, and on the negative $x_n$-axis when $b>0$, the spherical cap $B_b$ contains a hemisphere if $b<0$ and is contained in a hemisphere if $b>0$; for $b=0$, it is exactly a hemisphere.

\medskip

\emph{Throughout this section we assume that $b<0$ and set $\beta:=-b$.}

The following statement was proved in \cite{CPP} but only for $b\geq0$. We now extend it to $b<0$.

\begin{proposition}\label{compactness}
Assume that $b<0$ and  
\begin{equation}\label{strictmu1}
\mu_1(M)<\mu_1(B_b).
\end{equation}
Let $(u_i)$ be a sequence in $H^1(M)$ such that
\begin{equation}\label{PSweakzero}
\begin{cases}
u_i\rightharpoonup0\ \text{ weakly but not strongly in }H^1(M),\\
J(u_i)\to c<\mu_1(M)+\mu_1(B_b)\qquad\text{and}\qquad
J'(u_i)\to0 \text{ \ in \ }H^{-1}(M).
\end{cases}
\end{equation}
Then, after passing to a subsequence,
\begin{equation}\label{compactconclusions}
c\geq\mu_1(B_b)\qquad\text{and}\qquad
\dist_Q(u_i,\PP^+\cup\PP^-)\to0.
\end{equation}
\end{proposition}

\begin{proof}
We divide the proof into six steps.

\smallskip
\noindent\textbf{Step 1: Boundedness and nonvanishing of the interior critical mass.}
Since $\beta:=-b>0$, using the Palais--Smale hypotheses in \eqref{PSweakzero}, we see that, for some $C>0$,
\[ C(1+\|u_i\|_Q)\geq J(u_i)-\frac1pJ'(u_i)[u_i]=\Big(\frac{1}{2}-\frac{1}{p}\Big)Q(u_i)+\beta\Big(\frac{1}{q}-\frac{1}{p}\Big)\idm|u_i|^qd\sigma_g\geq\frac1n\|u_i\|_Q^2.
\] 
Therefore, $(u_i)$ is bounded in $H^1(M)$ and, consequently, $J'(u_i)[u_i]\to 0$. It follows that
\begin{equation}\label{NehariPS}
\|u_i\|_Q^2\leq Q(u_i)+\beta\int_{\partial M}|u_i|^q\,d\sigma_g
=a_n\int_M|u_i|^p\,dv_g+o(1).
\end{equation}
Thus, since $(u_i)$ does not converge strongly to $0$, there exists $\kappa_0>0$ such that a subsequence satisfies
\begin{equation}\label{internalmass}
\int_M|u_i|^p\,dv_g\geq\kappa_0\qquad\text{for all \ }i\in\n.
\end{equation}
Note also that, since $u_i\rightharpoonup0$ weakly in $H^1(M)$, due to the compactness of the Sobolev embedding and of the trace map, by passing to a subsequence we obtain
\begin{equation}\label{L2zero}
u_i\to0\quad\text{strongly in }L^2(M),\qquad
u_i|_{\partial M}\to0\quad\text{strongly in }L^2(\partial M).
\end{equation}

\smallskip
\noindent\textbf{Step 2: A local small-mass estimate.}
Choose $0<\kappa<\kappa_0$ small enough so that
\begin{equation}\label{kappasmall}
\max\{a_nS_M^{-1},C_*,C_\#\}\Big(\kappa^{2/n}+\kappa^{1/(n-1)}\Big)<\frac12,
\end{equation}
where $C_*$ and $C_\#$ are the positive constants defined in \eqref{eq:Cstar} and the proof of Lemma \ref{lem:A2} respectively, and
\[
S_M:=\inf_{v\in H^1(M)\setminus\{0\}}
\frac{Q(v)}{\left(\int_M|v|^p\,dv_g\right)^{2/p}}>0.
\]
Let $\eta\in C^\infty(M)$. From the identity
\[
Q(\eta u_i)=\langle u_i,\eta^2u_i\rangle_Q+\int_M|\nabla_g\eta|^2u_i^2\,dv_g
\]
we derive 
\begin{equation}\label{localidentity}
Q(\eta u_i)
+\beta\int_{\partial M}\eta^2|u_i|^q\,d\sigma_g
=a_n\int_M\eta^2|u_i|^p\,dv_g
+\int_M|\nabla_g\eta|^2u_i^2\,dv_g+J'(u_i)[\eta^2 u_i].
\end{equation}
Since $(u_i)$ is bounded in $H^1(M)$ and $J'(u_i)\to 0$ in $H^{-1}(M)$, we have that $J'(u_i)[\eta^2 u_i]=o(1)$.
 On the other hand, using H\"older's inequality and the definition of $S_M$, we get
\begin{align}\label{local L^p-bound}
a_n\int_M\eta^2|u_i|^p\,dv_g
&\leq a_n\left(\int_{\operatorname{supp}(\eta)}|u_i|^p\,dv_g\right)^{(p-2)/p}
\left(\int_M|\eta u_i|^p\,dv_g\right)^{2/p}\nonumber\\
&\leq\frac{a_n}{S_M}
\left(\int_{\operatorname{supp}(\eta)}|u_i|^p\,dv_g\right)^{2/n}Q(\eta u_i).
\end{align}
Therefore, assuming
\begin{equation} \label{eq:kappa}
\int_{\supp(\eta)}|u_i|^p\,dv_g\leq\kappa\qquad\text{for every \ }i\in\n,
\end{equation}
we derive from \eqref{L2zero}, \eqref{kappasmall}, \eqref{localidentity} and \eqref{local L^p-bound} that
\begin{equation}\label{epsilonregularity}
\frac12Q(\eta u_i)
+\beta\int_{\partial M}\eta^2|u_i|^q\,d\sigma_g=o(1).
\end{equation}

\smallskip
\noindent\textbf{Step 3: The radius of concentration tends to zero.}
Define
\begin{equation}\label{concentrationfunction}
\mathcal Q_i(r):=
 \max_{x\in\overline M}\left(\int_{B_r(x)}|u_i|^p\,dv_g + \int_{B_r(x)\cap\partial M}|u_i|^q\,d\sigma_g\right),
\end{equation}
where $B_r(x)$ is the open ball in $\overline M$ of radius $r$ and center $x$. By \eqref{internalmass}, for $\kappa$ as chosen in Step 2, there exist $\tilde x_i\in\overline M$ and $r_i>0$ such that
\begin{equation}\label{normalization}
\kappa=\mathcal Q_i(r_i)=
 \left(\int_{B_{r_i}(\tilde x_i)}|u_i|^p\,dv_g + \int_{B_{r_i}(\tilde x_i)\cap\partial M}|u_i|^q\,d\sigma_g\right).
\end{equation}
We claim that
\begin{equation}\label{rito0}
r_i\to0.
\end{equation}
If not, a subsequence would satisfy $r_i\geq r_0>0$. Then, we cover $\overline M$ with finitely many balls $B_{r_0/2}(y_j)$, $j=1,\ldots,m$, and choose functions $\eta_j\in C^\infty(M)$ that take the value $1$ in $B_{r_0/2}(y_j)$ and are supported in $B_{r_0}(y_j)$. From \eqref{normalization} we get
\begin{equation}
\int_{B_{r_0}(y_j)}|u_i|^p\,dv_g
\le \int_{B_{r_i}(y_j)}|u_i|^p\,dv_g
+\int_{B_{r_i}(y_j)\cap\partial M}|u_i|^q\,d\sigma_g
\le \mathcal{Q}_i(r_i)=\kappa,
\end{equation}
for all $ i\in\n, \ j=1,\ldots,m,$ that is, $\eta_j$ satisfies \eqref{eq:kappa}. Hence, \eqref{epsilonregularity} holds and it yields
\[
S_M\Big(\int_M|\eta_ju_i|^p\,dv_g\Big)^{2/p} \leq Q(\eta_ju_i)\to0\quad\text{as \ }i\to\infty\quad\text{for each \ } j=1,\ldots,m.
\]
Since $1\leq\sum_{j=1}^m|\eta_j|^p$ on $\overline M$, this implies that
$$\int_M|u_i|^p\,dv_g\to 0,$$
which contradicts \eqref{internalmass}. Thus, \eqref{rito0} holds.

\smallskip
\noindent\textbf{Step 4: Rescaling and nonvanishing of the limit profile.}
After passing to a subsequence,
\[
r_i^{-1}\dist_g(\tilde x_i,\partial M)\to d\in[0,\infty].
\]
We replace the points of concentration as follows. If $d<\infty$ we choose $x_i\in\partial M$ such that $\dist_g(\tilde x_i,x_i)=\dist_g(\tilde x_i,\partial M)$ and set $C_0:=d+2$, and if $d=\infty$ we set $x_i:=\tilde x_i$ and $C_0:=1$. Then, we have two cases:
\begin{itemize}
\item[$(1)$] $x_i\in\partial M$ for all $i\in \n$,
\item[$(2)$] $x_i\in M$ and $r_i^{-1}\dist_g(x_i,\partial M)\to \infty$.
\end{itemize}
In both cases $B_{r_i}(\tilde x_i)\subset B_{C_0r_i}(x_i)$. Therefore,
\begin{align} \label{eq:C_0}
&\int_{B_{C_0r_i}(x_i)}|u_i|^p\,dv_g + \int_{B_{C_0r_i}(x_i)\cap\partial M}|u_i|^q\,d\sigma_g \nonumber\\
&\qquad\qquad\geq\int_{B_{r_i}(\tilde x_i)}|u_i|^p\,dv_g + \int_{B_{r_i}(\tilde x_i)\cap\partial M}|u_i|^q\,d\sigma_g=\kappa.
\end{align}
We use normal coordinates at $x_i$ in case $(2)$ and  Fermi coordinates in case $(1)$. We fix $R_i>0$  such that these coordinates are defined in the ball  $B_{2R_i}(0)$ of radius $2R_i$ centered at $0$ in $\h$, where $\h$ is the half-space $\hn:=\{(z_1,\ldots,z_n)\in\rn:z_n\geq 0\}$ in case $(1)$ and it is $\rn$ in case $(2)$ and, abusing notation, in both cases we write $\exp_{x_i}$ for the coordinate chart. 
Additionally, we require that $r_i^{-1}R_i\to\infty$ and, in case $(2)$, that $R_i^{-1}\dist_g(x_i,\partial M)\to \infty$.

We fix a cut-off function $\zeta\in C^\infty_c(\rn)$ such that $0\leq\zeta\leq 1$, $\zeta(z)=1$ if $|z|<1$  and $\zeta(z)=0$ if $|z|>2$, and define
\begin{equation}\label{rescaling}
\widetilde u_i(z):=r_i^{\frac{n-2}{2}} u_i\bigl(\exp_{x_i}(r_i z)\bigr)\qquad\text{and}\qquad \zeta_i(z):=\zeta(R_i^{-1}r_iz).
\end{equation}
The pullback metric $g_i:=\exp_{x_i}^*g$ satisfies $C_1^{-1}\,dx\leq dv_{g_i}\leq C_1\,dx$ for some positive constant $C_1$. Then, using the rescaling $x=r_iz$  and recalling that $r_i^{-1}R_i\to \infty$, we obtain,
for $i$ large enough,
\begin{align*}
\int_{B_{C_0r_i}(x_i)}|u_i|^p\,dv_g\leq C_1\int_{B_{C_0r_i}(0)}\big|\zeta(R_i^{-1}x)\,u_i(\exp_{x_i}(x))\big|^p\,dx\leq C_1\int_{B_{2C_0}(0)}|\zeta_i\tilde u_i|^p\,dz.
\end{align*}
Similarly, in case $(1)$ we also have
\begin{align*}
\int_{B_{C_0r_i}(x_i)\cap\partial M}|u_i|^q\,d\sigma_g\leq C_1\int_{B_{2C_0}(0)\cap\partial\h}|\zeta_i\tilde u_i|^q\,ds.
\end{align*}
These inequalities, together with \eqref{eq:C_0}, yield, for $i$ large enough,
\begin{equation}\label{eq:C_0H}
0<\kappa\leq C_1\left(\int_{B_{2C_0}(0)}|\zeta_i\tilde u_i|^p\,dz + \int_{B_{2C_0}(0)\cap\partial\h}|\zeta_i\tilde u_i|^q\,ds\right).
\end{equation}
As $(u_i)$ is bounded in $H^1(M)$, a similar argument shows that $(\zeta_i\tilde u_i)$ is bounded in $D^{1,2}(\h)$. Passing to a subsequence, we have that 
\begin{equation}\label{weakprofile}
\zeta_i\widetilde u_i\rightharpoonup U\quad\text{weakly in \ }D^{1,2}(\h),\qquad \zeta_i\widetilde u_i\to U\quad\text{in \ }L^2_{loc}(\h)\qquad\zeta_i\widetilde u_i\to U\quad\text{a.e. in \ }\h.
\end{equation}
Next, we show that $U\neq 0$. 

Arguing by contradiction, assume that $U=0$. Fix $\vr\in(0,C_0)$ such that
$$\exp_{x_i}B_{r_i\vr}(r_iz)\subset B_{r_i}(\exp_{x_i}(r_iz))\qquad\text{for every \ }z\in \overline{B_{C_0}(0)}.$$
Let $z\in \overline{B_{2C_0}(0)}$ and $\eta\in C^\infty_c(B_\vr(z))$. Set $\eta_i(x):=\eta(r_i^{-1}\exp_{x_i}^{-1}(x))$. Then, for $i$ large enough, using \eqref{localidentity} we obtain
\begin{align}\label{eq:Unotnull1}
\int_\h|\nabla(\eta\zeta_i\tilde u_i)|^2\,dz &=\int_\h|\nabla(\eta\tilde u_i)|^2\,dz=\int_\h|\nabla((\eta_iu_i)\circ\exp_{x_i})|^2\,dx \nonumber\\
&\leq C_1\int_M|\nabla_g(\eta_iu_i)|^2\,dv_g \leq C_2\,Q(\eta_iu_i)\nonumber\\
&\leq C_2\Big(a_n\int_M\eta_i^2|u_i|^p\,dv_g
+\int_M|\nabla_g\eta_i|^2u_i^2\,dv_g+J'(u_i)[\eta^{2}_i u_i]\Big).
\end{align}
Now, from our choice of $\vr$, the Hölder and Sobolev inequalities, \eqref{concentrationfunction} and \eqref{normalization} we obtain
\begin{align}\label{eq:Unotnull2}
\int_M\eta_i^2|u_i|^p\,dv_g&\leq C_1\int_\h|\tilde u_i|^{p-2}(\eta\tilde u_i)^2\,dz\leq C_1\left(\int_{B_\vr(z)}|\tilde u_i|^p\,dz\right)^{2/n}\left(\int_\h|\eta\zeta_i\tilde u_i|^p\,dz\right)^{2/p}\nonumber\\
&\leq C_3S^{-1}\left(\int_{B_{r_i}(\exp_{x_i}(r_iz))}|u_i|^p\,dv_g\right)^{2/n}\int_\mathbb{H}|\nabla(\eta\zeta_i\tilde u_i)|^2\,dz\nonumber\\
&\leq C_3S^{-1}\kappa^{2/n}\int_\mathbb{H}|\nabla(\eta\zeta_i\tilde u_i)|^2\,dz.
\end{align}
Since $\zeta_i\widetilde u_i\to 0$ in $L^2_{loc}(\h)$, after rescaling we get that
\begin{align}\label{eq:Unotnull3}
\int_M|\nabla_g\eta_i|^2u_i^2\,dv_g\leq C_1\int_\h|\nabla\eta|^2(\zeta_i\tilde 
u_i)^2\,dz=o(1)
\end{align}
and $\|\eta^{2}_iu_i\|_{H^1(M)}$ is uniformly bounded. As $J'(u_i)\to 0$ in $H^{-1}(M)$, 
we have that $J'(u_i)[\eta^{2}_iu_i]=o(1)$. Then, setting $C_*:=C_2a_nC_3S^{-1}$, from 
\eqref{eq:Unotnull1}, \eqref{eq:Unotnull2} and \eqref{eq:Unotnull3} and our 
choice \eqref{kappasmall} of $\kappa$ we obtain
\begin{equation}\label{eq:Cstar}
\int_\mathbb{H}|\nabla(\eta\zeta_i\tilde u_i)|^2\,dz\leq C_*\kappa^{2/n}
\int_\mathbb{H}|\nabla(\eta\zeta_i\tilde u_i)|^2\,dz + o(1)<\frac{1}
{2}\int_\mathbb{H}|\nabla(\eta\zeta_i\tilde u_i)|^2\,dz + o(1).
\end{equation}
As a consequence, 
\[
\int_\mathbb{H}|\nabla(\eta\zeta_i\tilde u_i)|^2\,dz=o(1),\;\text{ yielding }\;
\int_\mathbb{H}|\eta\zeta_i\tilde u_i|^p\,dz=o(1),\; \int_{\partial\mathbb{H}}|\eta\zeta_i\tilde u_i|^q\,ds=o(1),
\]
 for every $\eta\in C^\infty_c(B_\vr(z))$ with $z\in \overline{B_{2C_0}(0)}$. Covering 
$\overline{B_{C_0}(0)}$ with a finite number of balls $B_{\vr/2}(z_1),\ldots,B_{\vr/2}
(z_m)$  with $z_j\in \overline{B_{2C_0}(0)}$, and choosing $\eta_j\in 
C^\infty_c(B_\vr(z_j))$ such that $\eta_j\equiv 1$ in $B_{\vr/2}(z_j)$ we conclude 
that 
\[
   %\begin{split}
   \int_{B_{2C_0}(0)\cap\mathbb H}|\zeta_i\widetilde u_i|^p\,dz
   +\int_{B_{2C_0}(0)\cap\partial\mathbb H}
     |\zeta_i\widetilde u_i|^q\,ds
\le
   \sum_{j=1}^m\left(
   \int_{\mathbb H}|\eta_j\zeta_i\widetilde u_i|^p\,dz
   +\int_{\partial\mathbb H}|\eta_j\zeta_i\widetilde u_i|^q\,ds
   \right)=o(1),
 %  \end{split}
   \]
contradicting \eqref{eq:C_0H}. This shows that $U\neq 0$.

\smallskip
\noindent\textbf{Step 5: The limiting equation.}
Let $J_\h:D^{1,2}(\h)\to\r$ be the functional 
$$
J_\h(w):=\frac{1}{2}\int_\h|\nabla w|^2\,dz-\frac{a_n}{p}\int_\h|w|^p\,dz+\frac{\beta}{q}\int_{\partial\h}| w|^q\,dz,
$$
(where the last integral is $0$ if $\h=\rn$). If $\vp\in C^\infty_c(\rn)$, set $\vp_i(\exp_{x_i}(z))=r_i^{(2-n)/2}\vp(r_i^{-1}z)$. Then, using \eqref{L2zero}, one sees that
$$
o(1)=J'(u_i)\vp_i=J_\h'(\zeta_i\tilde u_i)\vp +o(1),
$$
and passing to the limit we get that $J_\h'(U)=0$. That is, $U$ solves
\begin{equation}\label{halfspaceproblem}
\begin{cases}
-\Delta w=a_n|w|^{p-2}w&\text{in }\hn,\\
\frac{\partial w}{\partial z_n}=\beta|w|^{q-2}w&\text{on }\partial\hn
\end{cases}
\end{equation}
in case $(1)$, and it solves equation
\begin{equation}\label{wholeproblem}
-\Delta w=a_n|w|^{p-2}w\quad\text{in }\rn,
\end{equation}
in case $(2)$. Note that, since $\rn_+$ is conformally equivalent to the spherical cap $B_b$ via stereographic projection, we have 
$$J_{\rn_+}(U) \geq \mu_1(B_b),$$
and, similarly,
$$J_{\rn}(U) \geq \mu_1(\sn).$$

\smallskip
\noindent\textbf{Step 6: Subtraction of the limit profile.}
Let $\zeta$ be the cut-off function introduced in \eqref{rescaling} and define $v_i\in H^1(M)$ by
\[ 
v_i(x):=\begin{cases}
\zeta\big(R_i^{-1}\exp_{x_i}^{-1}(x)\big)\,r_i^\frac{2-n}{2}\,U\big(r_i^{-1}\exp_{x_i}^{-1}(x)\big) & \text{if \ }x\in\exp_{x_i}\big(B_{2R_i}(0)\big),
\\
0 & \text{otherwise},
\end{cases}
\]
and
$$w_i:=u_i-v_i.$$
Proposition \ref{prop:A} in the Appendix states that
\begin{equation}\label{profilesplit}
J(u_i)=J_\h(U)+J(w_i)+o(1)\qquad\text{and}\qquad J'(w_i)[w_i]=o(1).
\end{equation}
By \eqref{coerciveidentity}, applied to $w_i$,
\begin{equation*}
J(w_i)
=\frac1nQ(w_i)
+\beta\frac{n-2}{2n(n-1)}\int_{\partial M}|w_i|^q\,d\sigma_g+o(1)
\geq-o(1).
\end{equation*}
Therefore,  in case $(1)$, since $U$ is a nontrivial solution of \eqref{halfspaceproblem}, we have that
\begin{equation*}
J_\h(U)\geq\mu_1(B_b)
\end{equation*}
and,  in case $(2)$, since $U$ is a nontrivial solution of \eqref{wholeproblem}, we get
\begin{equation*}
J_\h(U)\geq\mu_1(\sn)>\mu_1(B_b),
\end{equation*}
where the last strict inequality is proven in \cite[Lemma A.1]{HanLi}.
Thus, in both cases, using assumptions \eqref{strictmu1} and \eqref{PSweakzero}, we obtain
\begin{equation}\label{eq:c>muB}
\mu_1(B_b)\leq J_\h(U)\leq c+o(1)<\mu_1(M)+\mu_1(B_b)<2\mu_1(B_b).
\end{equation}
It follows that $U$ does not change sign. Therefore,
$$v_i\in\PP^+\cup\PP^-.$$
Next, we show that $(w_i)$ converges strongly to $0$ in $H^1(M)$. Suppose that it does not. Then, by passing to a subsequence, we have that $w_i$ is bounded away from $0$ in $H^1(M)$. Setting $t_i:=t_{w_i}$, it follows from \eqref{eq:t_u} that $t_iw_i\in\N$. Given that $J'(w_i)[w_i]=o(1)$, we have $t_i\to1$. Therefore,
\begin{equation*}\label{remainderenergy}
J(w_i)\geq\mu_1(M)+o(1).
\end{equation*}
From this, we deduce that 
$$c\geq\mu_1(B_b)+\mu_1(M),$$
which contradicts assumption \eqref{PSweakzero}. Since $u_i-v_i=w_i$, we conclude that 
$$\dist_Q(u_i,\PP^+\cup\PP^-)\leq\|u_i-v_i\|_Q\to 0.$$
This, together with \eqref{eq:c>muB}, proves the assertions in \eqref{compactconclusions} and completes the proof of Proposition \ref{compactness}.
\end{proof}

\section{A criterion for attaining the least-energy nodal level}

For $b \geq 0$ the following abstract criterion for attaining $\mu_2(M)$ was proved in \cite[Theorem 2.8]{CPP}. Using Propositions \ref{PSconstruction} and \ref{compactness}, we can now extend it to $b < 0$. The proof is essentially the same; we include it for the sake of completeness.

\begin{theorem}\label{abstract}
Let $b<0$ and assume \eqref{strictmu1}. If
\begin{equation}\label{threshold}
\mu_2(M)<\mu_1(M)+\mu_1(B_b),
\end{equation}
then $\mu_2(M)$ is attained by a least-energy sign-changing solution of \eqref{P}.
\end{theorem}

\begin{proof}
Let $(u_k)$ be a sequence (as given by Proposition \ref{PSconstruction}) which satisfies
\begin{equation*}
J(u_k)\to\mu_2(M),\quad
\|J'(u_k)\|_{H^{-1}(M)}\to0,\quad\dist_Q(u_k,\E)\to0,\quad
\dist_Q(u_k,\PP^+\cup\PP^-)\geq\alpha>0.
\end{equation*}
The inequality
\begin{align*}
C(1+\|u_k\|_Q)\geq J(u_k)-\frac1qJ'(u_k)[u_k]
={}&\left(\frac12-\frac1q\right)Q(u_{k})
+a_n\left(\frac1q-\frac1p\right)\int_M|u_k|^p\,dv_g
\\
&\geq \frac{1}{2(n-1)}\|u_k\|_Q^2
\end{align*}
shows that $(u_k)$ is bounded. Passing to a subsequence,
\[
u_k\rightharpoonup u\quad\text{weakly in }H^1(M).
\]
Standard arguments, using this weak convergence and the Brezis--Lieb lemma, give
\begin{equation}\label{BL}
J'(u)=0,
\qquad
J(u_k)=J(u)+J(v_k)+o(1),
\qquad
J'(v_k)\to0,
\end{equation}
where $v_k:=u_k-u\rightharpoonup0$ weakly in $H^1(M)$. We claim that $(v_k)$ converges strongly in $H^1(M)$.

Suppose it does not. If $u\neq0$, then $u\in\N$ and $J(u)\geq\mu_1(M)$. Therefore
\[
J(v_k)\to\mu_2(M)-J(u)
\leq\mu_2(M)-\mu_1(M)<\mu_1(B_b),
\]
contradicting Proposition \ref{compactness}. On the other hand, if $u=0$, then $\dist_Q(v_k,\PP^+\cup\PP^-)\geq\alpha>0$. This, again, contradicts  Proposition \ref{compactness}. Hence, $u_k\to u$ strongly in $H^1(M)$.

Consequently, $J(u)=\mu_2(M)$, $J'(u)=0$ and, since $\dist_Q(u_k,\E)\to0$ and $\E$ is closed in $H^1(M)$, we have that $u\in\E$. Thus, $u$ is a least-energy nodal solution of \eqref{P}.
\end{proof}

Using the same argument employed in \cite[Corollary 2.9]{CPP}, we derive from Theorem \ref{abstract} a useful criterion for the existence of a least-energy sign-changing solution to \eqref{P}.

\begin{corollary}\label{pathcriterion}
Let $b<0$ and assume \eqref{strictmu1}. Assume further that there exist $u_0,u_1\in H^1(M)$, $u_0,u_1\geq0$, such that
\[
u_s:=(1-s)u_0-su_1\neq0\qquad\text{for every \ } s\in[0,1]
\]
and
\begin{equation}\label{pathbound}
J(ru_s)<\mu_1(M)+\mu_1(B_b)
\qquad\text{for every \ } r\geq0\text{ \ and \ }  s\in[0,1].
\end{equation}
Then \eqref{P} has a least-energy sign-changing solution.
\end{corollary}

\begin{proof}
For $R>0$, set $\sigma(s,t)=Rt u_s$. We claim that $\sigma \in \Sigma$. Indeed, it is immediate to verify that
$\sigma$ is continuous and that the boundary conditions in \eqref{Sigma} at
$s=0,\,1$ and $t=0$ are satisfied. In addition, as $u_s\neq0$,

\[
A_*:=\min_{s\in[0,1]}\int_M|u_s|^p\,dv_g>0.
\]
Thus, as $p>q>2$, we have that
\[
J(Ru_s)
\leq CR^2+C R^q-\frac{a_nA_*}{p}R^p\to-\infty,
\]
uniformly in $s$, and, as $\beta:=-b>0$, \ \eqref{tau} yields
\[
\tau(Ru_s)
=\frac{a_nR^p\int_M|u_s|^p}
{R^2Q(u_s)+\beta R^q\int_{\partial M}|u_s|^q}
\to\infty
\]
uniformly in $s$. 
Hence $R$ can be chosen so that $\sigma\in\Sigma$. Since the image of every $\sigma\in\Sigma$ intersects $\E$, using \eqref{pathbound} we obtain that
\[
\mu_2(M)<\mu_1(M)+\mu_1(B_b),
\]
and the statement follows from Theorem \ref{abstract}.
\end{proof}

\section{The test function for $b<0$}

Throughout this section, we assume that $b<0$, $n \geq 7$ and that there exists a nonumbilic point $\xi \in \partial M$. Our goal is to show that condition \eqref{threshold} holds under these hypotheses. To this end, we will use the test function introduced by Han and Li in \cite{HanLi}.

Following \cite{HanLi}, after a conformal change of metric, we may assume that
\begin{equation}\label{normalizationgeometry}
h_g(\xi)=0,\qquad R_{ij}(\xi)=0.
\end{equation}
We take Fermi coordinates at $\xi$ so that the second fundamental form of $\partial M$ at $\xi$ is diagonal, with principal curvatures
$\lambda_1,\dots,\lambda_{n-1}$. Since $h_g(\xi)=0$, \ 
$\sum_{j=1}^{n-1}\lambda_j=0$ 
and, since $\xi$ is nonumbilic, not all principal curvatures are $0$. Hence,
\begin{equation}\label{lambdanonzero}
\sum_{j=1}^{n-1}\lambda_j^2>0.
\end{equation}
Set
\begin{equation}\label{Tb}
T_b:=-\frac{b}{n-2}>0
\end{equation}
and consider the restriction of the bubble
\begin{equation}\label{bubble}
U_b(z):=
\bigl(1+|z'|^2+|z_n-T_b|^2\bigr)^{-\frac{n-2}{2}},
\qquad z=(z',z_n)\in\R^n_+,
\end{equation}
to the half-space $\rn_+:=\{(z_1,\ldots,z_n)\in\rn:z_n>0\}$. This is the positive solution to the boundary value problem
\begin{equation}\label{bubbleeq}
-\Delta w=a_n|w|^{p-2}w\quad\text{in }\R^n_+,
\qquad
\frac{\partial}{\partial z_n} w=-b|w|^{q-2}w\quad\text{on }\partial\R^n_+,
\end{equation}
whose energy functional is given by
$$
J_{\rn_+}(w):=\frac{1}{2}\ihn|\nabla w|^2-\frac{a_n}{p}\ihn|w|^p-\frac{b}{q}\idhn|w|^q.$$
Since the spherical cap $B_b$ is conformally equivalent to $\rn_+$ (via the stereographic projection), we have that $J(U_b)=\mu_1(B_b)$.

Let $\psi$ be a smooth cut-off function that takes the value $1$ in a neighborhood of the origin and is supported in $B_\rho(0)$, where $\rho>0$ is such that the Fermi coordinates at $\xi$ are defined in $B_{2\rho}(0)$, and, for $\eps>0$, define
\begin{equation}\label{Weps}
W_\eps(x):=\eps^{-\frac{n-2}{2}}
\psi(x)U_b(x/\eps).
\end{equation}
For $b<0$ and $n>5$, Han and Li proved that
\begin{equation*}
\max_{t>0}J(tW_\eps)
=\mu_1(B_b)+Q_3\eps^2+o(\eps^2).
\end{equation*}
In fact, their test function \cite[(1.8)]{HanLi}, which works for any $b\in\mathbb{R}$, contains an additional correction term $\delta\phi$, and they showed that this function satisfies
\[
\max_{t>0}J(tW_{\eps,\delta})
=\mu_1(B_b)+Q_1(\phi)\eps\delta+Q_2(\phi)\delta^2+Q_3\eps^2
+o(\eps^2+\delta^2).
\]
However, as they pointed out, for $b<0$ the correction is unnecessary because, in the notation of \cite[(2.9)]{HanLi},
\begin{equation*}
Q_3=-q_nT_b(1+T_b^2)^{\frac{3-n}{2}}
\sum_{j=1}^{n-1}\lambda_j^2,
\qquad
q_n:=\frac{\sigma_{n-2}(n-2)^2}{8(n-1)^2}
\int_0^\infty\frac{r^n}{(1+r^2)^{n-1}}\,dr>0
\end{equation*}
and, consequently, \eqref{Tb} and \eqref{lambdanonzero} imply that
\begin{equation}\label{Q3negative}
Q_3<0.
\end{equation}
Then, taking $\delta=0$ and $C_0:=-Q_3>0$, we obtain
\begin{equation}\label{strictbubble}
\max_{t>0}J(tW_\eps)
=\mu_1(B_b)-C_0\eps^2+o(\eps^2).
\end{equation}
This identity immediately yields $\mu_1(M)<\mu_1(B_b)$, which, as shown in \cite[Lemma 1.3]{HanLi}, implies the existence of a smooth positive solution $u_0$ to \eqref{P} satisfying
\begin{equation}\label{u0}
J(u_0)=\mu_1(M).
\end{equation}

Next, we will prove that $u_0$ and $u_1:=W_\eps$ satisfy the hypotheses of Corollary \ref{pathcriterion}. Specifically, we will first show that there exist $R>0$ and $\eps_1>0$ such that
\begin{equation}\label{eq:claim2}
J(su_0-tW_\eps)\leq 0 \quad\text{for all \ }s,t\in[0,r]\text{ \ with \ }s+t=r\geq R\text{ \ and \ }\eps\in (0,\eps_1).
\end{equation}
Afterwards, we will see that, for such an $R>0$, there exists an $\eps_0\in(0,\eps_1)$ such that
\begin{equation}\label{eq:claim1}
J(su_0-tW_\eps)<\mu_1(M)+\mu_1(B_b)\quad\text{for all \ }s,t\in[0,r]\text{ \ with \ }s+t=r\in[0,R]\text{ \ and \ }\eps\in(0,\eps_0).
\end{equation}

\subsection{The proof of statement \eqref{eq:claim2}}

Put
\[
v_{\theta,\eps}:=(1-\theta)u_0-\theta W_\eps,
\qquad0\leq\theta\leq1,\quad\eps>0.
\]
We claim that there are $A_0>0$ and $\eps_1>0$ such that
\begin{equation}\label{uniformAstar}
\int_M|v_{\theta,\eps}|^p\,dv_g\geq A_0
\qquad\text{for all \ }0<\eps<\eps_1\text{\ and \ }0\leq\theta\leq1.
\end{equation}
Suppose not. Then there are $\eps_j\to0$ and $\theta_j\to\theta\in[0,1]$ such that
$\|v_{\theta_j,\eps_j}\|_{L^p(M)}\to0$. On every compact subset of $\overline M\setminus\{\xi\}$ one has that $W_{\eps_j}\to0$ uniformly. Hence $(1-\theta)u_0=0$ there, and since $u_0>0$ we obtain $\theta=1$. On the other hand,
\[
\|W_\eps\|_{L^p(M)}^p\longrightarrow
\int_{\R^n_+}U_b^p\,dz>0.
\]
Thus, by the triangle inequality,
\[
\|v_{\theta_j,\eps_j}\|_{L^p(M)}
\geq \theta_j\|W_{\eps_j}\|_{L^p(M)}-(1-\theta_j)\|u_0\|_{L^p(M)},
\]
which stays bounded away from zero, a contradiction. This proves \eqref{uniformAstar}.

Note that $Q(W_\eps)$ and $\|W_\eps\|_{L^q(\partial M)}^q$ are uniformly bounded as $\eps\to0$. Consequently, there is $C>0$, independent of $\theta$ and small $\eps$, such that
\[
Q(v_{\theta,\eps})\leq C,
\qquad
\int_{\partial M}|v_{\theta,\eps}|^q\,d\sigma_g\leq C.
\]
Since $b<0$ and $p>q>2$, \eqref{uniformAstar} gives
\begin{align}
J(rv_{\theta,\eps})
&=\frac{r^2}{2}Q(v_{\theta,\eps})
-\frac{a_nr^p}{p}\int_M|v_{\theta,\eps}|^p\,dv_g
-\frac{b r^q}{q}\int_{\partial M}|v_{\theta,\eps}|^q\,d\sigma_g\notag\\
&\leq Cr^2+Cr^q-\frac{a_nA_0}{p}r^p.
\label{uniformlargeR}
\end{align}
Hence, there exists $R>0$ such that
\begin{equation}\label{largeRnegative}
J(rv_{\theta,\eps})\leq0
\qquad\text{for every \ }r\geq R,\quad 0\leq\theta\leq1\text{ \ and \ } 0<\eps<\eps_1,
\end{equation}
as stated in \eqref{eq:claim2}.

\subsection{The proof of statement \eqref{eq:claim1}}

We will use the following elementary fact: For any $\vr\geq2$ there exists $C_\vr>0$ such that
\begin{equation}\label{nonlinearineq}
\bigl||A-B|^\vr-A^\vr-B^\vr\bigr|
\leq C_\vr\bigl(A^{\vr-1}B+AB^{\vr-1}\bigr)
\qquad\text{for all \ }A,B\geq0.
\end{equation}
Fix $R>0$ and $\eps_1>0$ as above. For $0\leq s,t\leq R$, using \eqref{nonlinearineq} we get
\begin{align}
J(su_0-tW_\eps)
\leq{}&J(su_0)+J(tW_\eps)
+st\,|\langle u_0,W_\eps\rangle_Q|\notag\\
&+C_R\int_M\bigl(u_0^{p-1}W_\eps+u_0W_\eps^{p-1}\bigr)\,dv_g\notag\\
&+C_R\int_{\partial M}\bigl(u_0^{q-1}W_\eps+u_0W_\eps^{q-1}\bigr)\,d\sigma_g.
\label{masterinteraction}
\end{align}
Since $u_0$ is a positive solution of \eqref{P},
\begin{equation*}
\langle u_0,W_\eps\rangle_Q
=a_n\int_Mu_0^{p-1}W_\eps\,dv_g
+b\int_{\partial M}u_0^{q-1}W_\eps\,d\sigma_g,
\end{equation*}
so
\begin{equation}\label{crossbound}
|\langle u_0,W_\eps\rangle_Q|
\leq C\left(
\int_Mu_0^{p-1}W_\eps\,dv_g
+\int_{\partial M}u_0^{q-1}W_\eps\,d\sigma_g
\right).
\end{equation}

Near $\xi$, the bubble satisfies
\begin{equation*}
W_\eps(x)\leq
C\left(\frac{\eps}{\eps^2+|x|^2}\right)^{\frac{n-2}{2}}.
\end{equation*}
with $C$ depending on the fixed $b$. Since $u_0$ is smooth and positive, the first interior interaction satisfies
\begin{align*}
\int_Mu_0^{p-1}W_\eps\,dv_g
&\leq C\int_{B_\rho(0)}
\left(\frac{\eps}{\eps^2+|x|^2}\right)^{\frac{n-2}{2}}\,dx\leq C\eps^{\frac{n-2}{2}}.
\end{align*}
Indeed, after the change of variable $x=\eps y$, the integral over $|y|<\rho/\eps$ grows at most quadratically. Next, as $p-1=\frac{n+2}{n-2}$, we get
\begin{align*}
\int_Mu_0W_\eps^{p-1}\,dv_g
&\leq C\int_{\R^n}
\left(\frac{\eps}{\eps^2+|x|^2}\right)^{\frac{n+2}{2}}\,dx\leq C\eps^{\frac{n-2}{2}}.
\end{align*}
On the boundary, Fermi coordinates give
\begin{align*}
\int_{\partial M}u_0^{q-1}W_\eps\,d\sigma_g
&\leq C\eps^{\frac{n-2}{2}},\\ 
\int_{\partial M}u_0W_\eps^{q-1}\,d\sigma_g
&\leq C\int_{\R^{n-1}}
\left(\frac{\eps}{\eps^2+|x'|^2}\right)^{\frac n2}\,dx'
\leq C\eps^{\frac{n-2}{2}},
\end{align*}
where we used $q-1=\frac n{n-2}$. Consequently,
\begin{equation}\label{interactiontotal}
\int_M\bigl(u_0^{p-1}W_\eps+u_0W_\eps^{p-1}\bigr)\,dv_g
+\int_{\partial M}\bigl(u_0^{q-1}W_\eps+u_0W_\eps^{q-1}\bigr)\,d\sigma_g=O\big(\eps^{\frac{n-2}{2}}\big),
\end{equation}
and
\begin{equation}\label{interactionorder}
O\big(\eps^{\frac{n-2}{2}}\big)=o(\eps^2)
\quad\text{if }n\geq7.
\end{equation}
It is at this point that the restriction on the dimension of the manifold is required, because now, using \eqref{u0}, \eqref{strictbubble}, \eqref{crossbound}, \eqref{interactiontotal} and \eqref{interactionorder} to estimate the right-hand side of \eqref{masterinteraction}, we obtain
$$J(su_0-tW_\eps)\leq\mu_1(M)+\mu_1(B_b)-C_0\eps^2+o(\eps^2).$$
This proves the statement \eqref{eq:claim1}.

\subsection{The proof of Theorem \ref{main}}
Let $n\geq7$ and suppose that $\partial M$ has a nonumbilic point. Under these hypotheses, Han and Li proved \eqref{strictmu1}.
Assume further that $b<0$. Clearly,
$$u_s:=(1-s)u_0-sW_\eps\neq0\qquad\text{for all } s\in[0,1]\text{ and }\eps>0\text{ sufficiently small}.$$
Moreover, \eqref{eq:claim2} and \eqref{eq:claim1} imply that
\begin{equation*}
J(ru_s)<\mu_1(M)+\mu_1(B_b)
\qquad\text{for all } r\geq0 \text{ and } s\in[0,1].
\end{equation*}
Consequently, Corollary \ref{pathcriterion} establishes the existence of a least-energy sign-changing solution to \eqref{P}, as claimed by Theorem \ref{main}.
\medskip

\begin{remark}[The dimensions $n=5,6$]
\emph{The expansion by Han and Li used earlier is valid for $n>5$; thus, the case $n=5$ is not covered by the present approach. For $n=6$, our estimate of the nodal interaction is $O(\eps^{(n-2)/2}) = O(\eps^2)$, which is of the same order as the correction obtained by Han and Li; consequently, an analysis at the level of the coefficients, or a cancellation, would be required. Therefore, the present argument does not resolve the cases $n=5$ and $n=6$.}
\end{remark}

\appendix

\section{The energy splitting}
\label{appendix}

In this appendix we give the proof of the identities \eqref{profilesplit}. The data $u_i, r_i, x_i, \tilde u_i, U,$ etc. are those introduced in Section \ref{sec:compactness} and have the properties stated there. Recall that
\[ 
v_i(x):=\begin{cases}
\zeta\big(R_i^{-1}\exp_{x_i}^{-1}(x)\big)\,r_i^\frac{2-n}{2}\,U\big(r_i^{-1}\exp_{x_i}^{-1}(x)\big) & \text{if \ }x\in\exp_{x_i}\big(B_{2R_i}(0)\big),
\\
0 & \text{otherwise},
\end{cases}
\]
with $\zeta$ as in \eqref{rescaling}, and \ $w_i:=u_i-v_i.$ \ 
We begin with the following lemmas.

\begin{lemma}\label{lem:A1}
$v_i\rh 0\quad\text{weakly in \ }H^1(M)$.
\end{lemma}

\begin{proof}
Let $\vp\in C^\infty(M)$ and $r>0$. Then,
\begin{align*}
\Big|\int_{B_{r_ir}(x_i)}v_i\vp\,dv_g\Big|&\leq C_1\int_{B_{r_ir}(0)}\zeta(R_i^{-1}y)r_i^\frac{2-n}{2}|U(r_i^{-1}y)|\vp(\exp_{x_i}(y))\,dy \\
&\leq C_1|\vp|_\infty r_i^\frac{n+2}{2}\int_{B_{r}(0)}|U(z)|\,dz= o(1)\\ 
\Big|\int_{M\smallsetminus B_{r_ir}(x_i)}v_i\vp\,dv_g\Big|&\leq C_1|\vp|_\infty r_i^\frac{n+2}{2}\int_{B_{2R_ir_i^{-1}}(0)\smallsetminus B_{r}(0)}|U(z)|\,dz \\ 
&\leq C|\vp|_\infty R_i^\frac{n+2}{2} \Big(\int_\h|U(z)|^p\,dz\Big)^{1/p}=o(1).
\end{align*}
It follows that $\int_M v_i\vp\,dv_g\to 0$. Similarly, $\int_M \nabla_gv_i\cdot\nabla_g\vp\,dv_g\to 0$.
\end{proof}

\begin{lemma}\label{lem:A2}
For any $r>0$,  
$$\lim_{i\to\infty}\int_{B_r(0)}|\nabla(\tilde u_i-U)|^2=0.$$
\end{lemma}

\begin{proof}
Let $\vr>0$ be such that
$$\exp_{x_i}B_{r_i\vr}(r_iz)\subset B_{r_i}(\exp_{x_i}(r_iz))\qquad\text{for every \ }z\in \overline{B_{r}(0)},$$ 
and let $\eta\in C_c^\infty(B_\vr(z))$ for some $z\in \overline{B_r(0)}$. Then, using \eqref{weakprofile}, the mean value theorem, H\"older's and Sobolev's inequalities and the trace inequality, we see that
\begin{align*}
&\int_\h|\nabla(\eta(\tilde u_i-U))|^2= \int_\h\nabla(\tilde u_i-U)\cdot\nabla (\eta^2(\tilde u_i-U)) + o(1)\\
&=a_n\int_\h\big(|\tilde u_i|^{p-2}\tilde u_i-|U|^{p-2}U\big) (\eta^2(\tilde u_i-U))+b\int_{\partial\h}\big(|\tilde u_i|^{q-2}\tilde u_i-|U|^{q-2}U\big)(\eta^2(\tilde u_i-U)) + o(1)\\
&\leq C\Big(\int_\h\big(|\tilde u_i|^{p-2}+|U|^{p-2}\big) (\eta(\tilde u_i-U))^2+\int_{\partial\h}\big(|\tilde u_i|^{q-2}+|U|^{q-2}\big) (\eta(\tilde u_i-U))^2\Big)\\
&\leq C\left[\Big(\int_{B_\vr(z)}|\tilde u_i|^p\Big)^{2/n} + \Big(\int_{B_\vr(z)}|U|^p\Big)^{2/n}\right]\Big(\int_\h|\eta(\tilde u_i-U)|^p\Big)^{2/p}\\
&\qquad + C\left[\Big(\int_{B_\vr(z)\cap\partial\h}|\tilde u_i|^q\Big)^{1/(n-1)} + \Big(\int_{B_\vr(z)\cap\partial\h  }|U|^q\Big)^{1/(n-1)}\right]\Big(\int_{\partial\h}|\eta(\tilde u_i-U)|^q\Big)^{2/q}\\
&\leq C_\#\left[\Big(\int_{B_{r_i}(\exp_{x_i}(r_iz))}|u_i|^p\, dv_g\Big)^{2/n}+\Big(\int_{B_{r_i}(\exp_{x_i}(r_iz))\cap\partial M}|u_i|^q\, d\sigma_g\Big)^{1/(n-1)} + \eps_\vr\right]\int_\h|\nabla(\eta(\tilde u_i-U))|^2,
\end{align*}
where $\eps_\vr\to 0$ as $\vr\to 0$. Then, choosing $\vr$ sufficiently small, from \eqref{concentrationfunction}, \eqref{normalization}  and our choice \eqref{kappasmall} of $\kappa$ we derive
$$\int_\h|\nabla(\eta(\tilde u_i-U))|^2=o(1)\qquad\text{for every \ }\eta\in C_c^\infty(B_\vr(z))\text{ \ with \ }z\in \overline{B_{r}(0)}.$$
By reasoning as at the end of Step 4 of the proof of Proposition \ref{compactness}, we obtain the claim of the lemma.
\end{proof}

\begin{lemma}\label{lem:A3}
\begin{align*}
\int_M|\nabla_g w_i|^2\,dv_g &=\int_M|\nabla_g u_i|^2\,dv_g-\int_\h|\nabla U|^2\,dz + \theta(r,i) + o(1),
\\
\int_M|w_i|^p\,dv_g &=\int_M|u_i|^p\,dv_g-\int_\h|U|^p\,dz + \theta(r,i) + o(1),
\\
\int_{\partial M}|w_i|^q\,d\sigma_g &=\int_{\partial M}|u_i|^q\,d\sigma_g-\int_{\partial\h}|U|^q\,ds + \theta(r,i) + o(1),
\end{align*}
where
\begin{equation}\label{eq:theta}
\lim_{r\to\infty}\limsup_{i\to\infty}|\theta(r,i)|=0.
\end{equation} 
\end{lemma}

\begin{proof}
Once we have established Lemma \ref{lem:A2}, the proof follows the argument in \cite[pp. 40-41]{dhr}. We give the details for the sake of completeness. 

For any given $r>0$,
\begin{align*}
&\int_M|\nabla_g w_i|^2\,dv_g
\\
&\qquad=\int_{B_{r_ir}(x_i)}|\nabla_g w_i|^2\,dv_g+\int_{B_{2R_i}(x_i)\smallsetminus B_{r_ir}(x_i)}|\nabla_g w_i|^2\,dv_g+\int_{M\smallsetminus B_{2R_i}(x_i)}|\nabla_g u_i|^2\,dv_g.
\end{align*}
Set $\tilde g_i(z):=(\exp_{x_i}^*g)(r_iz)$. Then, using Lemma \ref{lem:A2}, we see that
\begin{equation}
\int_{B_{r_ir}(x_i)}|\nabla_g w_i|^2\,dv_g=\int_{B_r(0)}|\nabla(\tilde u_i-U)|^2\, dv_{\tilde g_i} = o(1).
\end{equation}
On the other hand, noting that
$$\lim_{r\to\infty}\limsup_{i\to\infty}\int_{B_{2R_i}(x_i)\smallsetminus B_{r_ir}(x_i)}|\nabla_g v_i|^2\,dv_g = 0,$$
and recalling that $(u_i)$ is bounded in $H^1(M)$, we see that
$$\int_{B_{2R_i}(x_i)\smallsetminus B_{r_ir}(x_i)}|\nabla_g w_i|^2\,dv_g=\int_{B_{2R_i}(x_i)\smallsetminus B_{r_ir}(x_i)}|\nabla_g u_i|^2\,dv_g + \theta(r,i),$$
where $\theta(r,i)$ satisfies \eqref{eq:theta}. Therefore,
$$\int_M|\nabla_g w_i|^2\,dv_g=\int_M|\nabla_g u_i|^2\,dv_g-\int_{B_{r_ir}(x_i)}|\nabla_g u_i|^2\,dv_g+ \theta(r,i)+o(1).$$
Since we are using normal coordinates at $x_i$ in case $(2)$ and Fermi coordinates in case $(1)$ we have that $\tilde g_i$ converges to the standard metric in $C^1(B_r(0))$. Therefore,
$$
\int_{B_{r_ir}(x_i)}|\nabla_g u_i|^2\,dv_g=
\int_{B_r(0)}|\nabla\tilde u_i|^2\, dv_{\tilde g_i}
= \int_{B_r(0)}|\nabla U|^2\,+o(1)
=\int_\h|\nabla U|^2\,dz+\eps_r+o(1),
$$
where $\eps_r\to0$ as $r\to\infty$. This shows that
$$\int_M|\nabla_g w_i|^2\,dv_g=\int_M|\nabla_g u_i|^2\,dv_g-\int_\h|\nabla U|^2\,dz + \theta(r,i) + o(1),$$
with $\theta(r,i)$ satisfying \eqref{eq:theta}. 

The other two identities are proved in a similar way.
\end{proof}

\begin{proposition}\label{prop:A}
\begin{equation*}
J(u_i)=J_\h(U)+J(w_i)+o(1)\qquad\text{and}\qquad J'(w_i)[w_i]=o(1),
\end{equation*}
\end{proposition}

\begin{proof}
By assumption (see Section \ref{sec:compactness}), $u_i\rh 0$ weakly in $H^1(M)$. Therefore, after passing to a subsequence,
\begin{equation*}
u_i\to0\quad\text{strongly in }L^2(M),\qquad\text{and}\qquad
u_i|_{\partial M}\to0\quad\text{strongly in }L^2(\partial M).
\end{equation*}
It follows that 
\begin{align*}
J(u_i)&=\frac12\int_M|\nabla_g u_i|^2\,dv_g-\frac{a_n}{p}\int_M|u_i|^p\,dv_g-\frac{b}{q}\int_{\partial M}|u_i|^q\,d\sigma_g + o(1), \\
J'(u_i)[u_i]&=\int_M|\nabla_g u_i|^2\,dv_g-a_n\int_M|u_i|^p\,dv_g-b\int_{\partial M}|u_i|^q\,d\sigma_g + o(1).
\end{align*}
By Lemma \ref{lem:A1}, the same is true for $w_i$. The statements now follow immediately from Lemma \ref{lem:A3}.
\end{proof}

\end{document}